\documentclass[11pt]{article}
\usepackage{graphicx}
\usepackage{amsmath}
\allowdisplaybreaks
\usepackage{amssymb}
\usepackage{amsthm}
\usepackage{epsfig}
\usepackage{mathrsfs}
\usepackage{amsbsy}
\usepackage{multirow}
\usepackage[bottom,symbol]{footmisc}

\usepackage[usenames,dvipsnames]{xcolor}
\usepackage[ruled, vlined, linesnumbered]{algorithm2e}
\let\oldnl\nl
\newcommand{\nonl}{\renewcommand{\nl}{\let\nl\oldnl}}
\usepackage{makecell}
\usepackage{changepage}
\usepackage{IEEEtrantools}
\usepackage{commath}
\usepackage{upgreek}
\usepackage{multicol}
\usepackage{subcaption}
\usepackage[margin=2.54cm]{geometry}
\usepackage{lipsum}
\usepackage{bigints}
\usepackage{framed}
\usepackage{colortbl}
\usepackage{calligra}
\usepackage{enumerate}
\usepackage{hyperref}
\usepackage{enumitem}
\usepackage{natbib}
\usepackage{abstract}
\setcitestyle{numbers,square,comma}
\usepackage{cleveref}
\usepackage{mathtools}
\usepackage{bbm}
\usepackage{appendix}
\usepackage{soul}
\usepackage{tikz}
\usepackage{pgfplots}
\usepackage{graphicx}
\usetikzlibrary{positioning}
\usepackage{array, multirow}

\usepackage{booktabs}
\newcolumntype{L}[1]{>{\raggedright\let\newline\\\arraybackslash\hspace{0pt}}m{#1}}
\newcolumntype{C}[1]{>{\centering\let\newline\\\arraybackslash\hspace{0pt}}m{#1}}
\newcolumntype{R}[1]{>{\raggedleft\let\newline\\\arraybackslash\hspace{0pt}}m{#1}}
\allowdisplaybreaks

\usepackage[margin = 1.0cm, font = small, labelfont = bf, labelsep = period]{caption}
\usepackage[font = small, labelfont = bf, labelsep = period]{caption}

\DeclareMathAlphabet{\mathcalligra}{T1}{calligra}{m}{n}

\newtheorem{theorem}{Theorem}

\newtheorem{corollary}{Corollary}
\newtheorem{lemma}{Lemma}

\newtheorem{proposition}{Proposition}

\newtheorem{definition}{Definition}

\newtheorem{remark}{Remark}

\newcommand{\bm}{\boldsymbol}
\newcommand{\bs}{\boldsymbol}

\newcommand{\C}{\tfrac{\pi^{p/2}}{\Gamma(p/2+1)}}

\newcommand{\mc}[1]{\mathcal{#1}}
\newcommand{\mb}[1]{\mathbb{#1}}
\newcommand{\mbb}[1]{\mathbbm{#1}}

\newcommand{\xit}{\widetilde{\bm{\xi}}}

\newcommand{\x}{{\bm{x}}}

\newcommand{\VV}{{\mathbb{V}}}

\newcommand{\EE}{{\mathbb{E}}}

\newcommand\redst{\bgroup\markoverwith{\textcolor{red}{\rule[0.5ex]{2pt}{0.4pt}}}\ULon}

\title{
Sample Complexity of Data-driven Multistage Stochastic Programming under Markovian Uncertainty
}

\author{ Hyuk Park and Grani A. Hanasusanto\\ Department of Industrial and Enterprise Systems Engineering\\ University of Illinois Urbana-Champaign, United States}

\date{}

\pgfplotsset{compat=1.18}

\begin{document}

\maketitle
\begin{abstract}
This work is motivated by the challenges of applying the sample average approximation (SAA) method to multistage stochastic programming with an unknown continuous-state Markov process. While SAA is widely used in static and two-stage stochastic optimization, it becomes computationally intractable in general multistage settings as the time horizon $T$ increases. Indeed, the number of samples required to obtain a reasonably accurate solution grows exponentially—a phenomenon known as the curse of dimensionality with respect to the time horizon. To overcome this limitation, we propose a novel data-driven approach, the Markov Recombining Scenario Tree (MRST) method, which constructs an approximate problem using only two independent trajectories of historical data. Our analysis demonstrates that the MRST method achieves polynomial sample complexity in $T$, providing a more efficient alternative to SAA. Numerical experiments on the Linear Quadratic Gaussian problem show that MRST outperforms SAA, addressing the curse of dimensionality.
\end{abstract}

\section{Introduction}
We consider the  multistage stochastic programming (MSP) problem with a time horizon $T$ given by
\begin{equation}\label{eq:true_nested}
    \min_{\substack{\bm x_1 \in \mc X_1(\bm x_0, \bm \xi_1) }} c(\bm x_1),
\end{equation}
\begin{equation}
\begin{aligned}
\text{where }c(\bm x_1):=c_1(\bm x_1,\bm\xi_1) + \mathbb{E} & \Bigg[ 
        \min _{ \substack{\bm x_2 \in \mathcal{X}_2(\bm x_1, \widetilde{\bm\xi}_2) }}  c_2(\bm x_2,\widetilde{\bm\xi}_2)
        \\
        &+\mathbb{E} \Bigg[
                    \cdots +
        \mathbb{E} \Bigg[ \min_{ \substack{ \bm x_T \in \mathcal{X}_T(\bm x_{T-1}, \widetilde{\bm \xi}_T)} } c_T(\bm x_T, \widetilde{\bm\xi}_T) \;\Bigg\vert\; \widetilde{\bm \xi}_{T-1} \Bigg] \cdots \;\Bigg\vert\; \widetilde{\bm \xi}_2 
                    \Bigg]  \;\Bigg\vert\; \bm \xi_1 
    \Bigg]
\end{aligned}
\label{eq:true_nested_2}
\end{equation}
is referred to as the nested formulation of the MSP problem \eqref{eq:true_nested}.
Here, $\bm x_t \in \mathbb{R}^{d_t}$ is the decision vector at stage $t$ for $t = 1, \ldots, T$, while  $\bm x_0 \in \mathbb{R}^{d_0}$ is given as a deterministic vector. The random vector $\widetilde{\bm \xi}_{t} \in \mathbb{R}^p$, supported on a set $\Xi\subseteq \mathbb{R}^p$, represents uncertain data revealed after the $t-1$ stage decision for $t = 2, \ldots, T$, while the first-stage data $\bm \xi_1$ is known, i.e., deterministic. Thus, $\widetilde{\bm \xi}_{[1:T]} := (\bm \xi_1, \widetilde{\bm \xi}_2, \ldots, \widetilde{\bm \xi}_T)$ constitutes a random data process. For $t = 1, \ldots, T$, $c_t: \mathbb{R}^{d_t} \times \mathbb{R}^p \rightarrow \mathbb{R}$ is a continuous function, and $\mc X_t: \mathbb{R}^{d_{t-1}} \times \mathbb{R}^p \rightrightarrows \mathbb{R}^{d_t}$ is a measurable multifunction. The first-stage feasible region $\mc X_1(\bm x_0, \bm \xi_1)$ is deterministic, given that $\bm x_0$ and $\bm \xi_1$ are deterministic.

The MSP problem can be viewed as a sequential decision-making process. At each time step $t = 1, \ldots, T$, the decision $\bm{x}_t$ is made based on the information available up to time $t$, i.e., ${\bm{\xi}}_{[1:t]} = (\bm{\xi}_1, \ldots, {\bm{\xi}}_t)$, while the future realizations $\widetilde{\bm{\xi}}_{[t+1:T]} = (\widetilde{\bm{\xi}}_{t+1}, \ldots, \widetilde{\bm{\xi}}_T)$ remain unknown. To solve the MSP problem \eqref{eq:true_nested}, one can adopt the dynamic programming equations equivalent to the nested formulation \eqref{eq:true_nested_2}, where subproblems are defined and solved in a backward manner for all $\bm{x}_{t-1}$ and $\bm{\xi}_t$. If the MSP problem is solved exactly, the decision-maker obtains the optimal first-stage decision (referred to as the here-and-now decision), $\bm{x}_1^\star$, and optimal policies for subsequent stages (referred to as wait-and-see decisions), $\bm{x}_t^\star(\bm{\xi}_{[1:t]})$, which are functions of all realized data up to time $t$ for $t = 2, \ldots, T$.

In this paper, we assume that the data process follows a time-homogeneous continuous-state Markov process, i.e., the random data $\widetilde{\bm\xi}_{t+1}$ given the history $\bm\xi_{[1:t]}$ up to $t$ stage does not depend on $\bm\xi_{[1:t-1]}$. Consequently, as shown in \eqref{eq:true_nested_2}, the \emph{conditional} expectation of $c_{t+1}(\cdot, \widetilde{\bm \xi}_{t+1})$ given $\bm \xi_t$ should be evaluated at each stage $t$. However, when the (conditional) distribution is continuous, even in a two-stage stochastic program (i.e., $T=2$), merely computing the (conditional) expectation for a fixed decision becomes challenging~(\citet{hanasusanto2016comment}). In principle, one can apply the sample average approximation (SAA) method, where the data process in the MSP problem \eqref{eq:true_nested} is approximated by a scenario tree generated using conditional Monte Carlo sampling ~(\citet[Chapter 5.8]{shapiro2021lectures}). In this scheme, $N_1$ samples of the random vector $\widetilde{\bm\xi}_2$, i.e., $\widehat{\bm\xi}_2^i$ for $i=1, \ldots, N_1$, are generated conditional on $\bm\xi_1$. Then, conditional on each $\widehat{\bm\xi}_2^i$ for $i=1,\ldots,N_1$, $N_2$ samples $\widehat{\bm\xi}_3^{ij}$ for $j=1, \ldots, N_2$ are generated, and this process is repeated up to the terminal stage $T$. In this way, the scenario tree contains $\sum_{t=1}^{T-1} \prod_{i=1}^{t} N_i$ samples.

From both practical and theoretical perspectives, several challenges arise when applying the SAA method to the MSP problem \eqref{eq:true_nested}. While the first-stage solution $\widehat{\bm x}_1$ of the SAA problem can always serve as an implementable (suboptimal) solution to \eqref{eq:true_nested}, the SAA policies $\widehat{\bm x}_t(\bm\xi_{[1:t]})$ for $t=2,\ldots,T$ are only implementable for the samples used in the scenario tree. Therefore, if the realization of 
 the random data $\widetilde{\bm\xi}_t$ does not match with one of the samples in the scenario tree, the decision maker should generate a new scenario tree with a time horizon of $T-t$ and re-solve the corresponding SAA problem to obtain an implementable solution at each time $t$. This approach may not be viable in cases with limited computational resources, such as onboard
computers, or in scenarios where decisions need to be made frequently, as in real-time control problems.

More importantly, the SAA method assumes that the true distribution of the data process is known—in particular, that the conditional distribution of $\widetilde{\bm \xi}_{t+1}$ given $\bm \xi_t$ is available in our setting. However, in many real-world scenarios, this exact distribution is unknown. Even if the distribution is known, numerically representing MSP requires some kind of discretization. In this regard, sampling-based schemes, including SAA, yield \emph{approximate} solutions to the MSP problem \eqref{eq:true_nested} and thus may result in suboptimality. Therefore, whichever method is adopted, it is essential to analyze its sample complexity—that is, the number of samples required to ensure that the approximate solution achieves the desired accuracy with high probability.

While the literature shows that accurate SAA solutions can be efficiently obtained for single-stage and two-stage cases, the analysis in \citet{kleywegt2002sample,shapiro2006complexity,shapiro2005complexity} suggests that, for $N_1 = \cdots = N_{T-1} = N$, the suboptimality of the SAA method is $\widetilde{\mc O}(T^{\frac{1}{2}} N^{-\frac{1}{2T}})$. This corresponds to a sample complexity of $\widetilde{\mathcal{O}}(T^T \epsilon^{-2T})$ to achieve $\epsilon$-suboptimality, which grows exponentially with the time horizon $T$—commonly referred to as \emph{the curse of dimensionality with respect to the time horizon}. Indeed, the authors argue that 
\begin{quote}
\emph{“Multistage stochastic programs are generically computationally intractable already when medium-accuracy solutions are sought.”}
\end{quote}
More recently, \citet{jiang2021complexity} and \citet{reaiche2016note} provide a similar sample complexity analysis for the SAA method and raise an open question of whether a provable scheme can be developed for MSP that overcomes the curse of dimensionality.

In this work, we address the aforementioned challenges. Our contributions can be summarized as follows:
\begin{enumerate}
    \item 
    We propose a novel data-driven approach, the Markov Recombining Scenario Tree (MRST) method, to formulate an approximate problem to the MSP problem with an underlying continuous-state Markov process \eqref{eq:true_nested}. Unlike the SAA method, MRST can provide implementable policy in a computationally efficient manner. Moreover, compared to existing methods, MRST accommodates a broader class of MSP problems and is more sample-efficient, requiring only two independent historical trajectories of the data process.
    \item 
    Under standard regularity conditions, we derive a non-asymptotic suboptimality bound for the MRST method and establish its sample complexity as $\widetilde{\mathcal{O}}(T^{p+3} \epsilon^{-p-2})$, demonstrating a polynomial dependence on the time horizon $T$. 
    To the best of our knowledge, there is no existing works that provide non-asymptotic sample complexity results for data-driven MSP problem under Markovian uncertainty. More importantly, the mild dependence of our sample complexity on $T$ directly addresses the curse of dimensionality.
    \item 
    We conduct stylized numerical experiments on the Linear Quadratic Gaussian (LQG) control problem. We compare the suboptimality of our scheme with the SAA method across varying time horizons. The results illustrate that our method successfully overcomes the curse of dimensionality.
\end{enumerate}

\paragraph{Notation}
Bold font $\bm \gamma$ and regular font $\gamma$ represent a vector and a scalar, respectively. The tilde symbol is used to denote randomness (e.g., $\widetilde{\bm \xi}_t$) to differentiate from their realizations (e.g., $\bm \xi_t$). 
The data process from stage $t$ to $t'$ is denoted as $\widetilde{\bm \xi}_{[t:t']}$.
For a random variable $\widetilde{\gamma}$, $\EE[\,\widetilde{\gamma}\,]$ and $\VV[\,\widetilde{\gamma}\,]$ denote its expectation and variance, respectively. Similarly, $\EE_{\xi}[\,\widetilde{\gamma}\,]:=\EE[\,\widetilde{\gamma}\,\vert\,\widetilde{\xi}=\xi]$ and $\VV_{\xi}[\,\widetilde{\gamma}\,]:=\VV[\,\widetilde{\gamma}\,\vert\,\widetilde{\xi}=\xi\,]$ denote conditional expectation and variance of a random variable $\widetilde{\gamma}$ given a realization $\xi$ of the random variable $\widetilde\xi$, respectively, For any $T \in \mb Z_{++}$, we define $[T]$ as the index set $\{1,\dots,T\}$.  The indicator function of an event $\mathcal{E}$ is defined through $\mbb I(\mathcal{E}) = 1$ if $\mathcal{E}$ is true and $\mbb I(\mathcal{E}) = 0$ otherwise. 
We use the standard big O notation: $\mathcal{O}(\cdot)$. In addition, $\widetilde{\mathcal{O}}(\cdot)$ is used to suppress multiplicative terms with logarithmic dependence. 

\section{Problem Statement}\label{sec:problem_statement}
Using dynamic programming principles~(\citet{bellman1957dynamic}), the MSP problem \eqref{eq:true_nested} can be reformulated in terms of a sequence of value functions $Q_t: \mathbb{R}^{d_{t-1}} \times \mathbb{R}^{p} \rightarrow \mathbb{R}$ for each stage $t \in [T]$:
\begin{equation}
\label{eq:true_dp_first}
\begin{array}{ccl}
Q_1(\bm x_{0}, \bm \xi_1) = & \textup{min} & \displaystyle 
 c_1(\bm x_1,\bm\xi_1) +  \EE_{\bm \xi_1} \Big[ Q_{2}(\bm x_{1}, \widetilde{\bm \xi}_{2} ) \Big] \\
& \textup{s.t.} & \displaystyle \bm x_1 \in \mc X_1(\bm x_{0}, \bm \xi_1),
\end{array}
\end{equation}
and, for each subsequent stage $t\in [T]\setminus \{1\}$, $\bm x_{t-1}\in \mc X_{t-1}(\bm x_{t-1}, \bm \xi_t)$, and $\bm\xi_t\in\Xi$, 
\begin{equation}
\label{eq:true_dp}
\begin{array}{ccl}
Q_t(\bm x_{t-1}, \bm \xi_t) = & \textup{min} & \displaystyle 
 c_t(\bm x_t,\bm\xi_t) +  \EE_{\bm \xi_t} \Big[ Q_{t+1}(\bm x_{t}, \widetilde{\bm \xi}_{t+1} ) \Big] \\
& \textup{s.t.} & \displaystyle \bm x_t \in \mc X_t(\bm x_{t-1}, \bm \xi_t).
\end{array}
\end{equation}
Here, $$\EE_{\bm \xi_t} \Big[ Q_{t+1}(\bm x_{t}, \widetilde{\bm \xi}_{t+1} ) \Big] := \EE \Big[ Q_{t+1}(\bm x_{t}, \widetilde{\bm \xi}_{t+1} ) \Big| \widetilde{\bm \xi}_{t} = \bm \xi_t \Big]$$ represents the conditional expectation of the value function at stage $t+1$, given the realization $\bm \xi_t$ at stage $t$.
For simplicity, we assume that $Q_{T+1}(\bm x_{T}, \bm \xi_{T+1}) = 0$ for all $\bm x_{T}$ and $\bm \xi_{T+1}$, which implies that no additional cost is incurred beyond the terminal stage $T$.

\subsection{Markov Recombining Scenario Tree Method}
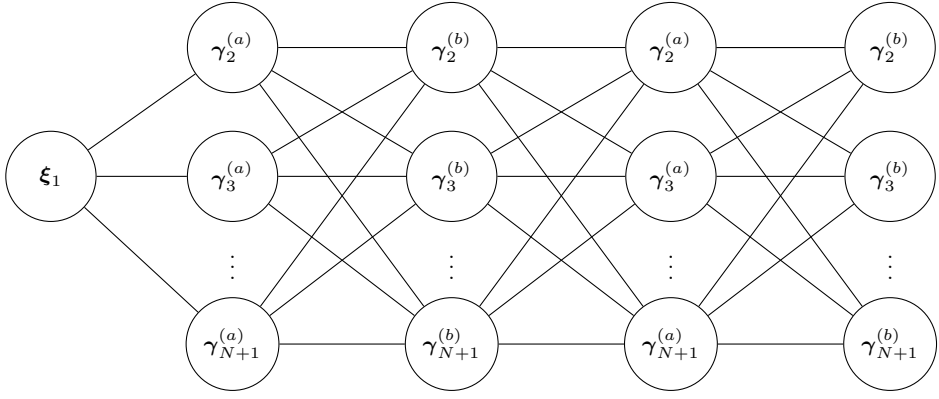
\begin{figure}
\centering
        \begin{tikzpicture} 
        [
            new_node/.style={circle, draw, minimum size=1.2cm}, 
            empty_node/.style={fill=white, minimum size=1.2cm}] 
            \node[new_node]  (t-2)                              {\scriptsize $\bm \gamma_3^{(a)}$};
            \node[new_node]  (t-1)    [above= 0.5cm of t-2]     {\scriptsize $\bm \gamma_2^{(a)}$};
            \node[empty_node]  (t-:)    [below= 0.1cm of t-2]     {\scriptsize };
            \node[new_node]  (t-N)    [below= 1.0cm of t-2]     {\scriptsize $\bm \gamma_{N+1}^{(a)}$};
            \node[new_node]    (3-2)    [left= 1.7cm of t-2]     {\scriptsize $\bm \gamma_3^{(b)}$};
            \node[new_node]    (3-1)    [above= 0.5cm of 3-2]     {\scriptsize $\bm \gamma_2^{(b)}$} ;
            \node[fill=white]  (3-:)    [below= 0.1cm of 3-2]     {\scriptsize $\vdots$} ;
            \node[new_node]    (3-N)    [below= 1.0cm of 3-2]     {\scriptsize $\bm \gamma_{N+1}^{(b)}$} ;   
            \node[fill=white]  (t-:)    [below= 0.1cm of t-2]     {\scriptsize $\vdots$} ;
            \node[new_node]    (2-2)    [left= 1.7cm of 3-2]     {\scriptsize $\bm \gamma_3^{(a)}$};
            \node[new_node]    (2-1)    [above= 0.5cm of 2-2]     {\scriptsize $\bm \gamma_2^{(a)}$} ;
            \node[fill=white]  (2-:)    [below= 0.1cm of 2-2]     {\scriptsize $\vdots$} ;
            \node[new_node]    (2-N)    [below= 1.0cm of 2-2]     {\scriptsize $\bm \gamma_{N+1}^{(a)}$} ;
            \node[new_node]    (1)      [left= 1.2cm of 2-2]      {\scriptsize $\bm \xi_1$};
            \node[new_node]    (T_1-2)    [right= 1.7cm of t-2]     {\scriptsize $\bm \gamma_3^{(b)}$};
            \node[new_node]    (T_1-1)    [above= 0.5cm of T_1-2]     {\scriptsize $\bm \gamma_2^{(b)}$} ;
            \node[fill=white]  (T_1-:)    [below= 0.1cm of T_1-2]     {\scriptsize $\vdots$} ;
            \node[new_node]    (T_1-N)    [below= 1.0cm of T_1-2]     {\scriptsize $\bm \gamma_{N+1}^{(b)}$} ;
            \path[-] (1) edge node {} (2-1);
            \path[-] (1) edge node {} (2-2);
            \path[-] (1) edge node {} (2-N);
            \path[-] (2-1) edge node {} (3-1);
            \path[-] (2-1) edge node {} (3-2);
            \path[-] (2-1) edge node {} (3-N);
            \path[-] (2-2) edge node {} (3-1);
            \path[-] (2-2) edge node {} (3-2);
            \path[-] (2-2) edge node {} (3-N);
            \path[-] (2-N) edge node {} (3-1);
            \path[-] (2-N) edge node {} (3-2);
            \path[-] (2-N) edge node {} (3-N);
            \path[-] (3-1) edge node {} (t-1);
            \path[-] (3-1) edge node {} (t-2);
            \path[-] (3-1) edge node {} (t-N);
            \path[-] (3-2) edge node {} (t-1);
            \path[-] (3-2) edge node {} (t-2);
            \path[-] (3-2) edge node {} (t-N);
            \path[-] (3-N) edge node {} (t-1);
            \path[-] (3-N) edge node {} (t-2);
            \path[-] (3-N) edge node {} (t-N);
            \path[-] (t-1) edge node {} (T_1-1);
            \path[-] (t-1) edge node {} (T_1-2);
            \path[-] (t-1) edge node {} (T_1-N);
            \path[-] (t-2) edge node {} (T_1-1);
            \path[-] (t-2) edge node {} (T_1-2);
            \path[-] (t-2) edge node {} (T_1-N);
            \path[-] (t-N) edge node {} (T_1-1);
            \path[-] (t-N) edge node {} (T_1-2);
            \path[-] (t-N) edge node {} (T_1-N);
            \path[-] (t-1) edge node {} (T_1-1);
        \end{tikzpicture}
        \vspace{0cm}
        \caption{Scenario tree using two historical trajectories when $T=5$.}
        \label{fig:recombining_tree}
    \end{figure}

We propose a method to approximate the unknown Markovian data process, which we call the Markov Recombining Scenario Tree (MRST). It is based on two independent realized trajectories—denoted as \(\bm{\gamma}^{(a)}_{[1:N+1]}\) and \(\bm{\gamma}^{(b)}_{[1:N+1]}\)—of the underlying Markov process \(\widetilde{\bm{\xi}}_{[1:N+1]}\). As depicted in \Cref{fig:recombining_tree}, MRST is constructed by alternating between the two trajectories over time: at even-numbered stages, nodes are drawn from \(\bm{\gamma}^{(a)}_{[2:N+1]}\) of the first trajectory, while at odd-numbered stages, nodes are taken from \(\bm{\gamma}^{(b)}_{[2:N+1]}\) of the second trajectory—the first elements for the trajectories, i.e., $\bm{\gamma}^{(a)}_1$ and $\bm{\gamma}^{(b)}_1$, are excluded.  We then approximate the conditional expectation \(\mathbb{E}_{\bm \xi_t}[Q_{t+1}(\bm{x}_t, \widetilde{\bm{\xi}}_{t+1})]\) via the following empirical estimates:
\noindent for all odd-numbered $t\in[T]$,
\begin{align} 
    &\widehat{\EE}_{\bm \xi_t} \Big[ 
    Q_{t+1}(\bm x_{t}, \widetilde{\bm \xi}_{t+1} ) 
    \Big] 
    =\begin{cases} 
    \displaystyle
    \sum_{i\in[N]} Q_{t+1}(\bm x_{t}, {\bm \gamma}_{i+1}^{(a)} )\cdot
     \left(\frac{ \mc K_h(\bm \xi_t - {\bm \gamma}_i^{(a)}) }{\sum_{j \in [N]} 
    \mc K_h(\bm \xi_t - {\bm \gamma}_j^{(a)}) 
    }\right)
    & \text { if } \displaystyle \sum_{j=1}^N \mc K_h(\bm \xi_t - {\bm \gamma}_j^{(a)}) >0 
        \\ 
        0 & \text { otherwise,}\end{cases}  \label{eq:odd_approx_expec}
\end{align}
and, for all even-numbered $t\in[T]$,
\begin{align} 
    &\widehat{\EE}_{\bm \xi_t} \Big[ 
    Q_{t+1}(\bm x_{t}, \widetilde{\bm \xi}_{t+1} )\Big] 
     =\begin{cases} 
    \displaystyle
    \sum_{i\in[N]} Q_{t+1}(\bm x_{t}, {\bm \gamma}_{i+1}^{(b)} )\cdot
     \left(\frac{ \mc K_h(\bm \xi_t - {\bm \gamma}_i^{(b)}) }{\sum_{j \in [N]} 
    \mc K_h(\bm \xi_t - {\bm \gamma}_j^{(b)}) 
    }\right)
    & \text { if } \displaystyle \sum_{j=1}^N \mc K_h(\bm \xi_t - {\bm \gamma}_j^{(b)})  >0 
        \\ 
        0 & \text { otherwise.}
        \end{cases}
    \label{eq:even_approx_expec}
\end{align}
Here, $\mc K_h:\mathbb{R}^{p}\times\mathbb{R}_{+}\rightarrow\mathbb{R}_{+}$ represents a kernel function with a parameter \( h \geq 0 \) known as the bandwidth parameter. In this paper, we consider the following
\begin{equation}
    \mc K_h(\bm \theta) = k(\Vert\bm \theta\Vert) \cdot \mb I(\Vert\bm \theta\Vert\leq h),
    \label{eq:kernel}
\end{equation}
where $\Vert\cdot\Vert$ denotes the Euclidean norm, $k:[0,h]\rightarrow\mathbb{R}_{+}$ is a non-increasing function, and $\mathbb{I}(\cdot)$ is an indicator function that returns 1 if $\Vert\bm \theta\Vert$ is less than or equal to $h$, and 0 otherwise. 

Note that, in \eqref{eq:odd_approx_expec} and \eqref{eq:even_approx_expec} we approximate the conditional expectation of the \emph{true} value functions \(Q_{t+1}(\cdot)\) in \eqref{eq:true_dp_first} and \eqref{eq:true_dp}, which are not available without precise knowledge of the underlying Markov process. Hence, we further need to define the \emph{approximate} value function \(\widehat{Q}_{t+1}(\cdot)\). By recursively applying approximate conditional expectations from the terminal stage backward to the initial stages, we define our approximate problem through the following dynamic programming equations:
\begin{equation}
\label{eq:approx_dp}
\begin{array}{ccl}
\widehat{Q}_1(\bm x_{0}, \bm \xi_1) = & \textup{min} & \displaystyle 
 c_1(\bm x_1,\bm\xi_1) +  \widehat{\EE}_{\bm \xi_1} \Big[ \widehat{Q}_{2}(\bm x_{1}, \widetilde{\bm \xi}_{2} ) \Big] \\
& \textup{s.t.} & \displaystyle \bm x_1 \in \mc X_1(\bm x_{0}, \bm \xi_1),
\end{array}
\end{equation}
for even-numbered $t\in[T]\setminus\{1\}$ and for all $i\in[N+1]\setminus\{1\}$
\begin{equation}
\label{eq:approx_dp_even}
\begin{array}{ccl}
\widehat{Q}_t(\bm x_{t-1}, \bm \gamma_{i}^{(a)}) = & \textup{min} & \displaystyle 
 c_t(\bm x_t,\bm \gamma_{i}^{(a)}) +  \widehat{\EE}_{\bm \gamma_{i}^{(a)}} \Big[ \widehat{Q}_{t+1}(\bm x_{t}, \widetilde{\bm \xi}_{t+1} ) \Big] \\
& \textup{s.t.} & \displaystyle \bm x_t \in \mc X_t(\bm x_{t-1}, \bm \gamma_{i}^{(a)}),
\end{array}
\end{equation}
and, for odd-numbered $t\in[T]\setminus\{1\}$ and for all $i\in[N+1]\setminus\{1\}$
\begin{equation}
\label{eq:approx_dp_odd}
\begin{array}{ccl}
\widehat{Q}_t(\bm x_{t-1}, \bm \gamma_{i}^{(b)}) = & \textup{min} & \displaystyle 
 c_t(\bm x_t,\bm \gamma_{i}^{(b)}) +  \widehat{\EE}_{\bm \gamma_{i}^{(b)} } \Big[ \widehat{Q}_{t+1}(\bm x_{t}, \widetilde{\bm \xi}_{t+1} )  \Big] \\
& \textup{s.t.} & \displaystyle \bm x_t \in \mc X_t(\bm x_{t-1}, \bm \gamma_{i}^{(b)}),
\end{array}
\end{equation}
 with \(\widehat{Q}_{T+1}(\bm{x}_T, \bm{\xi}_{T+1}) = 0\) for all \(\bm{x}_T\) and \(\bm{\xi}_{T+1}\). Let us refer to this dynamic program as the MRST problem. 

\begin{remark}
    Suppose that \(\widehat{\bm{x}}_t(\bm{\xi}_{[1:t]})\) for \(t = 2, \ldots, T\) are the optimal policies for the MRST problem. Note that these policies are defined only for the historical observations \(\bm{\gamma}^{(a)}_{[1:N+1]}\) and \(\bm{\gamma}^{(b)}_{[1:N+1]}\), i.e., in-sample data, as the true data process is approximated by our scenario tree (\Cref{fig:recombining_tree}). Consequently, these in-sample policies are not directly implementable in the true environment, where the data process may differ from the in-sample data. Nevertheless, our approach allows us to derive an implementable policy \(\widehat{\bm{x}}_t(\bm{\xi}_{[1:t]})\) for any unseen data. 
    Given the history \(\bm{\xi}_{[1:t]}\) and the approximate value function \(\widehat{Q}_{t+1}(\cdot)\),
    we can compute the out-of-sample policy \(\widehat{\bm{x}}_t(\bm{\xi}_{[1:t]})\) by first calculating the conditional expectation given $\bm \xi_t$ using the estimator in \eqref{eq:odd_approx_expec} or \eqref{eq:even_approx_expec}, and then solving the following minimization problem:
    \[
    \min_{\bm{x}_t} \; c_t(\bm{x}_t, \bm{\xi}_t) + \widehat{\mathbb{E}}_{\bm{\xi}_t } \Big[ \widehat{Q}_{t+1}(\bm{x}_t, \widetilde{\bm{\xi}}_{t+1}) \Big] \quad \text{s.t.} \quad \bm{x}_t \in \mathcal{X}_t(\bm{x}_{t-1}, \bm{\xi}_t).
    \]
    In contrast, the SAA method requires generating a new scenario tree with a \(T-t\) time horizon conditional on \(\bm{\xi}_t\) and solving the corresponding SAA problem each time new data is observed. This makes the SAA method not only reliant on the assumption of a known conditional distribution but also computationally inefficient compared to our method for deriving an implementable policy.
\end{remark}

\section{Sample Complexity Analysis}\label{sec:out-of-sample-performance}
We shift our focus to the theoretical analysis of the proposed approach. We formally state our assumptions about the MSP problem \eqref{eq:true_dp_first}.
\begin{enumerate}[label=(\textbf{A\arabic*})]
    \item \label{a3}  
    The first-stage feasible region $\mc X_1(\bm x_0, \bm \xi_1)$ is non-empty and has a finite diameter $D_1$. Similarly, for each subsequent stage $t \in [T] \setminus \{1\}$, the feasible region $\mc X_t(\bm x_{t-1}, \bm \xi_t)$ is non-empty and there exists a positive constant $D_t$ such that for every $\bm x_{t-1} \in \mc X_{t-1}(\cdot)$ and $\bm \xi_t \in \Xi$, $\mc X_t(\bm x_{t-1}, \bm \xi_t)$ has a finite diameter $D_t$. 
    \item \label{assumption:lipschitz_value}
    For $t\in[T]$, $Q_{t} (\x_{t-1}, \bm \xi_{t} )$ lies between 0 and 1. 
    Furthermore, there exists a Lipschitz constant $L_t$ such that
      \begin{equation*}
        \Big|Q_{t+1} (\x, \bm \xi_{t+1} )-Q_{t+1}(\bm y, \bm \xi_{t+1})\Big| \leq  L_{t+1} \Vert \x-\bm y \Vert, \quad \forall \x,\bm y \in \mathcal{X}_{t}\left(\bm x_{t-1}, \bm \xi_{t}\right) \; 
        \forall \bm \xi_{t+1}\in \Xi.
    \end{equation*}
    \item \label{a2} 
     $\widetilde{\bm{\xi}}_{[1:T]}$ is a time-homogeneous continuous-state Markov process supported on a bounded set $\Xi$.
      The conditional density function $f(\bm{\xi}|\bm{\gamma})$, defined for all $\bm\xi,\bm\gamma\in\Xi$, governs the Markov process. We assume that $f(\bm{\xi}|\bm{\gamma})$ is square-integrable and non-zero at any $\bm\gamma\in\Xi$, i.e., there exists $\underline{f}$ such that $0<\underline{f}\leq f(\bm{\xi}|\bm{\gamma})$ for all $\bm \gamma \in\Xi$.
    \item \label{a7_lip} 
    For $t\in[T]$, there exists a Lipschitz constant $M_{t+1}$ such that, for any $\bm\gamma,\bm\gamma'\in\Xi$
    \begin{align*}
        &\left| \EE_{\bm\gamma} \left[ Q_{t+1}(\bm x_{t}, \widetilde{\bm \xi}_{t+1})  \right] -  \EE_{\bm\gamma'} \left[ Q_{t+1}(\bm x_{t}, \widetilde{\bm \xi}_{t+1}) \right] \right|\leq M_{t+1}\Vert \bm\gamma - \bm\gamma' \Vert,\\
        &
        \left| \EE_{\bm\gamma} \left[ Q_{t+1}(\bm x_{t}, \widetilde{\bm \xi}_{t+1})^2  \right] -  \EE_{\bm\gamma'} \left[ Q_{t+1}(\bm x_{t}, \widetilde{\bm \xi}_{t+1})^2  \right] \right|\leq M_{t+1}\Vert \bm\gamma - \bm\gamma' \Vert.
    \end{align*}
    \item \label{assumption_kernel}
    The function $k(\cdot)$ in \eqref{eq:kernel} is non-increasing and bounded above by a constant $\overline{k}>0$, i.e., $k(\Vert\bm \theta\Vert) \leq \overline{k}$ for all $\bm \theta$.
\end{enumerate}
\noindent 
\ref{a3} and \ref{assumption:lipschitz_value} are standard assumptions in previous works analyzing the sample complexity of the MSP problem~(\citet{shapiro2006complexity,shapiro2005complexity}). These assumptions ensure the boundedness of the value functions. The assumption that the value functions lie between 0 and 1 is not restrictive—given their boundedness, they can always be scaled and translated to fall within the interval $[0, 1]$ by appropriately adjusting the cost functions $c_{t}(\cdot)$ for all $t \in [T]$ without altering the optimal policy.
A function $k(\cdot)$ in \eqref{eq:kernel} that satisfies \ref{assumption_kernel} includes several commonly used kernels with bounded support, such as the uniform kernel, $\mc K_h(\theta) = (1/h) \cdot \mathbb{I}(\Vert \bm \theta \Vert \leq h)$, and the Epanechnikov kernel, $\mc K_h(\theta) = (3/4)(h^2 - \Vert \bm \theta \Vert^2) \cdot \mathbb{I}(\Vert \bm \theta \Vert \leq h)$~(\citet{tibshirani2013nonparametric}).
The following definition of a Markov process in \citet[Proposition 7.6]{kallenberg1997foundations} will be useful later.
\begin{definition}
\label{def:markov}
    Let $\widetilde{\bm \xi}_{[1:\infty]}$ be a process with values in a Borel space $\Xi$. Then, $\widetilde{\bm \xi}_{[1:\infty]}$ is a time-homogeneous Markov process if and only if
    there exists some measurable function $\phi:\Xi\times[0,1]\rightarrow\Xi$ and i.i.d.~random variables $\zeta_t\sim Unif(0,1)$ independent of $\bm\xi_1$ such that $\bm\xi_{t+1}=\phi(\bm\xi_t,\zeta_{t+1})$ almost surely for all $t=1,2,\ldots.$.
\end{definition}

To highlight the main result, we first present the sample complexity of the MRST problem~\eqref{eq:approx_dp} in the following. 
\begin{theorem}[Sample Complexity of MRST]
\label{corol:sample_comp}
    Suppose that ${\bm x}_1^\star$ is an optimal solution to the MSP problem~\eqref{eq:true_dp_first}, and $\widehat{\bm x}_1^N$ is an optimal solution to the MRST problem~\eqref{eq:approx_dp}, constructed using two i.i.d.~Markovian trajectories $\widetilde{\bm \gamma}^{(a)}_{[1:N+1]}$ and $\widetilde{\bm \gamma}^{(b)}_{[1:N+1]}$, each of length $N+1$. Then, for any fixed $\epsilon\in[0,1]$, we have 
    \begin{equation*}
        c_{{1}}(\bm x_1^\star,\bm\xi_1)
        + \EE_{\bm \xi_1} \left[ Q_2 (\bm x_1^\star, \widetilde{\bm \xi}_{2} )  \right]
         \leq         
         c_{{1}}(\widehat{\bm x}_1,\bm\xi_1)
        + \EE_{\bm \xi_1} \left[ Q_{2}(\widehat{\bm x}_1,  \xit_{2}) \right] 
        +
        \epsilon
\end{equation*}
with high probability as long as 
$N \geq {\color{black}
\widetilde{\mathcal{O}}(T^{p+3} \epsilon^{-p-2})
}$.
\end{theorem}
\noindent As claimed in \Cref{corol:sample_comp}, the sample complexity of our approach exhibits only a polynomial dependence on $T$ and, therefore, does not suffer from the curse of dimensionality with respect to the time horizon. In what follows, we present the key derivations that establish \Cref{corol:sample_comp}. Specifically, we begin by deriving the generalization bound for the estimator in \eqref{eq:odd_approx_expec} and \eqref{eq:even_approx_expec}. Building on this, we extend the analysis to obtain the generalization and suboptimality bounds, which demonstrate the performance of the optimal policy for the MRST problem. Lastly, we compare the sample complexity of our approach with that of the SAA method.

The main challenge in our analysis arises from the Markovian structure of the unknown data process. In our data-driven setup, we merely have access to two trajectories of sequential data, where the samples within each trajectory are non-i.i.d.\ due to the underlying Markovian dependence. 
In contrast, existing works assume that the conditional distribution is known, allowing the generation of i.i.d.~conditional samples at each stage $t$ given any current realization ${\bm \xi}_{t}$~(\citet{jiang2021complexity,shapiro2006complexity}).
This simpler setup facilitates the direct application of tools from statistical learning theory. Unfortunately, their analyses reveal that the SAA method under such conditional sampling suffers from the curse of dimensionality with respect to the time horizon $T$, i.e., the number of samples required to achieve a fixed suboptimality $\epsilon$ grows exponentially with $T$. To the best of our knowledge, no existing work has proposed an alternative to the SAA method that effectively addresses this challenge with non-asymptotic analysis.
We begin by establishing preliminary results concerning the conditional expectation of the value functions.
\begin{lemma}\label{corol:expect_lipschitz}
For each stage $t\in[T]$, for any $\bm \xi_{t}$,
 there exists a constant $L_{t+1}\geq0$ such that
\begin{align*}
\left| \EE_{\bm \xi_{t}} \left[ Q_{t+1}(\x, \widetilde{\bm \xi}_{t+1}) \right] 
-  
\EE_{\bm \xi_{t}} \left[ Q_{t+1}(\bm y, \widetilde{\bm \xi}_{t+1}) \right] \right| & \leq  L_{t+1} \Vert \x-\bm y \Vert, \quad\quad\quad \forall \x,\bm y \in \mathcal{X}_{t}\left(\bm x_{t-1}, \bm \xi_{t}\right), 
\\
\left| \widehat\EE_{\bm \xi_{t}} \left[ Q_{t+1}(\x, \widetilde{\bm \xi}_{t+1}) \right] 
-
\widehat\EE_{\bm \xi_{t}} \left[ Q_{t+1}(\bm y, \widetilde{\bm \xi}_{t+1}) \right] \right| & \leq  L_{t+1} \Vert \x-\bm y \Vert, \quad\quad\quad \forall \x,\bm y \in \mathcal{X}_{t}\left(\bm x_{t-1}, \bm \xi_{t}\right),  \\
 \left| \EE_{\bm \xi_{t}} \left[\widehat Q_{t+1}(\x, \widetilde{\bm \xi}_{t+1})\right] 
 -  
 \EE_{\bm \xi_{t}} \left[\widehat Q_{t+1}(\bm y, \widetilde{\bm \xi}_{t+1}) \right] \right| & \leq  L_{t+1} \Vert \x-\bm y \Vert, \quad\quad\quad \forall \x,\bm y \in \mathcal{X}_{t}\left(\bm x_{t-1}, \bm \xi_{t}\right),
 \\
\left|\widehat \EE_{\bm \xi_{t}} \left[\widehat Q_{t+1}(\x, \widetilde{\bm \xi}_{t+1})\right] 
-  
\widehat\EE_{\bm \xi_{t}} \left[\widehat Q_{t+1}(\bm y, \widetilde{\bm \xi}_{t+1})\right] \right| & \leq  L_{t+1} \Vert \x-\bm y \Vert, \quad\quad\quad \forall \x,\bm y \in \mathcal{X}_{t}\left(\bm x_{t-1}, \bm \xi_{t}\right).
\end{align*}
\end{lemma}
\begin{proof}
    The Lipschitz continuity in \ref{a7_lip} directly implies the following: for each stage $t\in[T-1]$, we have
    \begin{align*}
        &\left| \EE_{\bm \xi_{t}} \left[ Q_{t+1}(\bm x, \widetilde{\bm \xi}_{t+1})  \right] -  \EE_{\bm \xi_{t}} \left[ Q_{t+1}(\bm y, \widetilde{\bm \xi}_{t+1})\right] \right| \\
        = & \bigg| \int_{\bm \xi_{t+1} \in \Xi}
        Q_{t+1}(\bm x, \bm \xi_{t+1}) \cdot f(\bm \xi_{t+1} | \bm \xi_t) \, \mathrm{d}\bm \xi_{t+1} 
        - 
        \int_{\bm \xi_{t+1} \in \Xi}
        Q_{t+1}(\bm y, \bm \xi_{t+1}) \cdot f(\bm \xi_{t+1} | \bm \xi_t) \, \mathrm{d}\bm \xi_{t+1} 
        \bigg|
        \\
        \leq &
         \int_{\bm \xi_{t+1} \in \Xi} \bigg|
        Q_{t+1}(\bm x, \bm \xi_{t+1}) 
        - 
        Q_{t+1}(\bm y, \bm \xi_{t+1}) 
        \bigg| \cdot f(\bm \xi_{t+1} | \bm \xi_t) 
        \, \mathrm{d}\bm \xi_{t+1} 
        \\
        \leq & L_{t+1} \Vert \bm x - \bm y \Vert.
    \end{align*}
    Furthermore, for all odd-numbered $t \in [T]$, using the definitions of the estimator in \eqref{eq:odd_approx_expec}, if $\sum_{j=1}^N \mathcal{K}_h(\bm \xi_t - {\bm \gamma}_j^{(a)}) > 0$, we have
    \begin{align*}
    &\left| \widehat{\EE}_{\bm \xi_{t}} \left[ Q_{t+1}(\bm x, \widetilde{\bm \xi}_{t+1}) \right] -  \widehat{\EE}_{\bm \xi_{t}} \left[ Q_{t+1}(\bm y, \widetilde{\bm \xi}_{t+1})\right] \right| \\
        = &  \Bigg|
        \sum_{i \in [N]} Q_{t+1}(\bm x, \bm \gamma_{i+1}^{(a)} ) \cdot
     \left(\frac{ \mathcal{K}_h(\bm \xi_t - \bm \gamma_i^{(a)}) }{\sum_{j \in [N]} 
    \mathcal{K}_h(\bm \xi_t - \bm \gamma_j^{(a)}) 
    }\right)
    -
    \sum_{i \in [N]} Q_{t+1}(\bm y, \bm \gamma_{i+1}^{(a)} ) \cdot
     \left(\frac{ \mathcal{K}_h(\bm \xi_t - \bm \gamma_i^{(a)}) }{\sum_{j \in [N]} 
    \mathcal{K}_h(\bm \xi_t - \bm \gamma_j^{(a)}) 
    }\right)
    \Bigg|
    \\
    \leq &
    \sum_{i \in [N]} \bigg| Q_{t+1}(\bm x, \bm \gamma_{i+1}^{(a)} )
    -
     Q_{t+1}(\bm y, \bm \gamma_{i+1}^{(a)} )\bigg|
     \cdot
     \left(\frac{ \mathcal{K}_h(\bm \xi_t - \bm \gamma_i^{(a)}) }{\sum_{j \in [N]} 
    \mathcal{K}_h(\bm \xi_t - \bm \gamma_j^{(a)}) 
    }\right)
    \\
    \leq & L_{t+1} \Vert \bm x - \bm y \Vert.
    \end{align*}
    Otherwise, if $ \sum_{j=1}^N \mathcal{K}_h(\bm \xi_t - \bm \gamma_j^{(a)}) = 0$ the above inequality holds trivially. Similarly, the remaining inequalities in \Cref{corol:expect_lipschitz} can be 
    verified.
\end{proof}

\noindent Next, we present results on the conditional variance of the value functions, stated in the following.

\begin{lemma}\label{a7}  
      For each stage $t \in [T]$, there exists a constant $\sigma_{t+1} > 0$ such that for any $\bm x_{t} \in \mc X_{t}(\bm x_{t-1}, \bm \xi_t )$ and $\bm \xi_t\in\Xi$, it holds that $\VV_{\bm \xi_{t}} [ Q_{t+1}(\bm x_t, \widetilde{\bm \xi}_{t+1}) ] \leq \sigma_{t+1}^2$.
\end{lemma}
\begin{proof}
    The boundedness condition in \ref{assumption:lipschitz_value} and the square-integrability in \ref{a2} ensure that the conditional variance $\VV_{\bm \xi_{t}} [ Q_{t+1}(\bm x_t, \widetilde{\bm \xi}_{t+1})]$ is finite for any $\bm x_t \in \mc X_t(\bm x_{t-1}, \bm \xi_t)$ and $\bm \xi_t \in \Xi$ at each stage $t \in [T]$.
\end{proof}
\begin{lemma}
\label{lem:lip_variance}
    For any $t \in [T]$ and $\bm\gamma, \bm\gamma' \in \Xi$, there exists a constant $M_{t+1}\geq0$ such that
    \begin{equation}
        \left| \VV_{\bm\gamma} \left[ Q_{t+1}(\bm x_{t}, \widetilde{\bm \xi}_{t+1})  \right] 
        -  
        \VV_{\bm\gamma'} \left[ Q_{t+1}(\bm x_{t}, \widetilde{\bm \xi}_{t+1}) \right] \right|\leq 3M_{t+1}\Vert \bm\gamma - \bm\gamma' \Vert.
        \label{eq:lem_variance}
    \end{equation}
\end{lemma}
\begin{proof} Using the definition of conditional variance, we can expand \eqref{eq:lem_variance} as follows: 
    \begin{align}
        &\bigg| \VV_{\bm\gamma} \left[ Q_{t+1}(\bm x_{t}, \widetilde{\bm \xi}_{t+1}) \right] -  \VV_{\bm\gamma'} \left[ Q_{t+1}(\bm x_{t}, \widetilde{\bm \xi}_{t+1}) \right] \Bigg|
        \nonumber
        \\ 
        = & 
        \left| 
        \left(\EE_{\bm\gamma} \left[ Q_{t+1}(\bm x_{t}, \widetilde{\bm \xi}_{t+1})^2 \right] 
        -  
        \EE_{\bm\gamma} \left[ Q_{t+1}(\bm x_{t}, \widetilde{\bm \xi}_{t+1}) \right]^2
        \right)
        -
        \left(\EE_{\bm\gamma'} \left[ Q_{t+1}(\bm x_{t}, \widetilde{\bm \xi}_{t+1})^2\right] -  
        \EE_{\bm\gamma'} \left[ Q_{t+1}(\bm x_{t}, \widetilde{\bm \xi}_{t+1})\right]^2
        \right)
        \right|
        \nonumber
        \\
        \leq &
        \left| 
        \EE_{\bm\gamma} \left[ Q_{t+1}(\bm x_{t}, \widetilde{\bm \xi}_{t+1})\right]^2
        -
        \EE_{\bm\gamma'} \left[ Q_{t+1}(\bm x_{t}, \widetilde{\bm \xi}_{t+1})\right]^2
        \right|
        +
        \left| 
        \EE_{\bm\gamma} \left[ Q_{t+1}(\bm x_{t}, \widetilde{\bm \xi}_{t+1})^2 \right] 
        -  
        \EE_{\bm\gamma'} \left[ Q_{t+1}(\bm x_{t}, \widetilde{\bm \xi}_{t+1})^2\right]
        \right|
        \nonumber
        \\
        \leq &
        \left| 
        \EE_{\bm\gamma} \left[ Q_{t+1}(\bm x_{t}, \widetilde{\bm \xi}_{t+1}) \right] 
        +  
        \EE_{\bm\gamma'} \left[ Q_{t+1}(\bm x_{t}, \widetilde{\bm \xi}_{t+1})\right]
        \right| 
        \cdot
        \left| 
        \EE_{\bm\gamma} \left[ Q_{t+1}(\bm x_{t}, \widetilde{\bm \xi}_{t+1}) \right] 
        -  
        \EE_{\bm\gamma'} \left[ Q_{t+1}(\bm x_{t}, \widetilde{\bm \xi}_{t+1})\right]
        \right| 
        \nonumber
        \\
        & + 
        \left| 
        \EE_{\bm\gamma} \left[ Q_{t+1}(\bm x_{t}, \widetilde{\bm \xi}_{t+1})^2\right] 
        -  
        \EE_{\bm\gamma'} \left[ Q_{t+1}(\bm x_{t}, \widetilde{\bm \xi}_{t+1})^2\right]
        \right|
        \nonumber
        \\
        \leq & 
        \left(
            \left| 
                \EE \left[ Q_{t+1}(\bm x_{t}, \widetilde{\bm \xi}_{t+1}) \middle\vert \bm\gamma \right] 
            \right|
            +
            \left| 
                \EE \left[ Q_{t+1}(\bm x_{t}, \widetilde{\bm \xi}_{t+1}) \middle\vert \bm\gamma' \right] 
            \right|
        \right)
        \cdot
        \left| 
        \EE \left[ Q_{t+1}(\bm x_{t}, \widetilde{\bm \xi}_{t+1}) \middle\vert \bm\gamma \right] 
        -  
        \EE \left[ Q_{t+1}(\bm x_{t}, \widetilde{\bm \xi}_{t+1}) \middle\vert \bm\gamma' \right]
        \right| \nonumber
        \\
        & +  
        \left| 
        \EE \left[ Q_{t+1}(\bm x_{t}, \widetilde{\bm \xi}_{t+1})^2 \middle\vert \bm\gamma \right] 
        -  
        \EE \left[ Q_{t+1}(\bm x_{t}, \widetilde{\bm \xi}_{t+1})^2 \middle\vert \bm\gamma' \right]
        \right|
        \nonumber
        \\
        \leq & 2\cdot M_{t+1}\Vert \bm \gamma - \bm \gamma'  \Vert
        + M_{t+1}\Vert \bm \gamma - \bm \gamma'  \Vert \label{eq:12}
        \\
        \leq & 3\cdot M_{t+1}\Vert \bm \gamma - \bm \gamma'  \Vert.
        \nonumber
    \end{align}
    Here, the inequality \eqref{eq:12} holds due to the boundedness condition in \ref{assumption:lipschitz_value} and the Lipschitz continuity in \ref{a7_lip}.
\end{proof}

\noindent As discussed at the beginning of this section, the main challenge lies in addressing the non-i.i.d.~samples consisting of a trajectory of observed data, e.g., $\bm \gamma^{(a)}_1$, $\ldots$, $\bm \gamma^{(a)}_{N+1}$ are non-i.i.d. The following proposition addresses this difficulty and presents a result on the generalization errors of the estimators \eqref{eq:odd_approx_expec} and \eqref{eq:even_approx_expec}, which serves as the key ingredient for our sample complexity analysis.

\begin{proposition}\label{prop:kernel_estimation_error}
    Suppose that $\bm \gamma\in\Xi$ is the current given state of a Markov process and $\widetilde{\bm \gamma}_{[1:N+1]}$ denotes a single Markovian trajectory where the initial state $\widetilde{\bm \gamma}_{1}$ follows the marginal distribution $\mu(\mathrm{d}\bm\gamma_1)$. Consider the estimator defined as follows: 
    \begin{equation}
    \label{eq:naive_E_hat}
        \widehat \EE_{\bm\gamma} \left[ Q_{t+1}(\bm x_{t}, \widetilde{\bm \xi}_{t+1})\right]
        = \begin{cases}\displaystyle\sum_{i=1}^N Q_{t+1}(\bm x_{t}, \widetilde{\bm \gamma}_{i+1}) \cdot
        \left(\frac{ \mc K_h( \bm \gamma - \widetilde{\bm \gamma}_{i}) }{\sum_{j \in [N]} 
    \mc K_h(\bm \gamma - \widetilde{\bm \gamma}_{j}) 
    }\right) & \text { if } \displaystyle \sum_{j \in [N]} 
    \mc K_h(\bm \gamma - \widetilde{\bm \gamma}_{j}) >0 
        \\ 
        0 & \text { otherwise. }\end{cases}
    \end{equation}
   Then, for any fixed $\epsilon\in[0,1]$, $h\geq0$, and $\bm x_t\in\mc X_t(\bm x_{t-1}, \bm \xi_t)$, we have
   \begin{equation}
       \left| \EE_{\bm\gamma} \left[ Q_{t+1}(\bm x_{t}, \widetilde{\bm \xi}_{t+1})\right] - \widehat \EE_{\bm\gamma} \left[ Q_{t+1}(\bm x_{t}, \widetilde{\bm \xi}_{t+1})\right] \right| \leq M_{t+1} h + \epsilon
   \end{equation}
  with probability at least 
   \begin{equation}
   1-2\exp\left(-\frac{NC_ph^p\underline{f}\underline{k}}{ 2\alpha}\cdot\frac{\epsilon^2}{\overline{k}\VV_{\bm\gamma}\left[ Q_{t+1}(\bm x_{t}, \widetilde{\bm \xi}_{t+1})\right]+3\overline{k}M_{t+1}h+\frac{1}{3}\epsilon}\right).
   \label{eq:prob_kernel_estimation_errors}
   \end{equation}
  Here, $C_p = \C$, and $\alpha > 0$ is a universal constant. Additionally, $\underline{f}>0$ and $M_{t+1}\geq 0$ are constants defined in \ref{a2} and \ref{a7_lip}, respectively, while $\overline{k}>0$ and $\underline{k}\geq0$ are constants that depend only on the choice of a function $k(\cdot)$ in \eqref{eq:kernel}.
\end{proposition}
\begin{proof}
    We first introduce the following function \
    \begin{equation}
    \label{eq:function_m}
        \widetilde{m}(\bm \gamma,\widetilde{\bm \gamma}_{[1:N]})
        = \begin{cases}\displaystyle\sum_{i=1}^N \EE_{\widetilde{\bm \gamma}_{i}}\left[ Q_{t+1}(\bm x_{t}, \widetilde{\bm \xi}_{t+1}) \right] \cdot
        \left(\frac{ \mc K_h(\bm \gamma - \widetilde{\bm \gamma}_{i}) }{\sum_{j \in [N]} 
    \mc K_h(\bm \gamma - \widetilde{\bm \gamma}_{j}) 
    }\right) & \text { if } \displaystyle \sum_{j \in [N]} 
    \mc K_h(\bm \gamma - \widetilde{\bm \gamma}_{j}) >0 
        \\ 
        0 
        & \text { otherwise, }\end{cases}
    \end{equation}
    where $\bm{\gamma}$ is the given current state (hence, no randomness), and $\widetilde{\bm{\gamma}}_{[1:N]}$ are the first $N$ components of a single random trajectory.
    Define a subset $\Xi_{h}(\bm \gamma)=\{\bm \xi\in\Xi:\Vert \bm\xi - \bm \gamma \Vert \leq h\}$ and its complement $\Xi^{^c}_h(\bm \gamma)=\Xi\setminus\Xi_{h}(\bm \gamma)$ where $\Xi_h(\bm \gamma)$ can be viewed as the intersection between the support set $\Xi$ and an Euclidean norm ball of radius $h$ centered at $\bm \gamma$.
    To ease the notation in the subsequent derivation, let us further define the following events: ${\mathcal{E}}_i=\{\widetilde{\bm \gamma}_{i} \in \Xi_h(\bm \gamma)\}$ and ${\mathcal{E}}_i^{^c}=\{\widetilde{\bm \gamma}_{i} \in \Xi_h^{^c}(\bm \gamma)\}$, $\forall i \in [N]$.
    Then, the indicator function $\mathbb{I}(\Vert \bm \gamma - \widetilde{\bm \gamma}_i\Vert\leq h )$ can be rewritten as $\mathbb{I}({\mathcal{E}}_i)$ for all  $i \in [N]$.  
    By triangle inequality, we have
    \begin{multline}
    \label{eq:decomposition}
        \left| \EE_{\bm \gamma} \left[ Q_{t+1}(\bm x_{t}, \widetilde{\bm \xi}_{t+1})  \right] - \widehat \EE_{\bm \gamma} \left[ Q_{t+1}(\bm x_{t}, \widetilde{\bm \xi}_{t+1}) \right] \right| 
        \\
        \leq
        \underbrace{
        \left| \EE_{\bm \gamma} \left[ Q_{t+1}(\bm x_{t}, \widetilde{\bm \xi}_{t+1})  \right] - \widetilde{m}(\bm \gamma,\widetilde{\bm \gamma}_{[1:N]}) \right| 
        }_{(a)}
        +
        \underbrace{
         \left| \widetilde{m}(\bm \gamma,\widetilde{\bm \gamma}_{[1:N]})
         -
         \widehat{\EE}_{\bm \gamma} \left[ Q_{t+1}(\bm x_{t}, \widetilde{\bm \xi}_{t+1})\right] \right| 
         }_{(b)}.
    \end{multline}
    Our goal is to derive upper bounds on the terms (a) and (b), separately and combine the results to complete the proof. 
    
    \noindent\underline{$(a)$:}
    Following the definition of \eqref{eq:kernel}, we have 
         \begin{align}
       (a)= & \Bigg\vert 
       \left( 
       {\EE}_{\bm \gamma} \left[ Q_{t+1}(\bm x_{t}, \widetilde{\bm \xi}_{t+1})\right]
       -
       \sum^N_{i=1}{\EE}_{\widetilde{ \bm \gamma}_{i}} \left[ Q_{t+1}(\bm x_{t}, \widetilde{\bm \xi}_{t+1})  \right]
       \frac{k(\Vert
       \bm \gamma-\widetilde{\bm \gamma}_{i}
       \Vert) \cdot \mb I({\mathcal{E}}_i)}{\sum^N_{j=1}k(\Vert\bm \gamma-\widetilde{\bm \gamma}_{j}
       \Vert) \cdot \mb I({\mathcal{E}}_j)}
       \right)\mathbb{I}(\cup^N_{i=1}{\mathcal{E}}_i)  
       \nonumber
       \\
       & + 
       \left(
       {\EE}_{\bm \gamma} \left[ Q_{t+1}(\bm x_{t}, \widetilde{\bm \xi}_{t+1}) \right]-0
       \right)\mathbb{I}(\cap^N_{i=1}{\mathcal{E}}^{^c}_i)
       \Bigg\vert  \label{eq:def_m}
       \\
       \leq &
       \Bigg\vert
       \left( 
       {\EE}_{\bm \gamma} \left[ Q_{t+1}(\bm x_{t}, \widetilde{\bm \xi}_{t+1})\right]
       -
       \sum^N_{i=1}{\EE}_{\widetilde{ \bm \gamma}_{i}} \left[ Q_{t+1}(\bm x_{t}, \widetilde{\bm \xi}_{t+1})  \right]
       \frac{k(\Vert
       \bm \gamma-\widetilde{\bm \gamma}_{i}
       \Vert) \cdot \mb I({\mathcal{E}}_i)}{\sum^N_{j=1}k(\Vert\bm \gamma-\widetilde{\bm \gamma}_{j}
       \Vert) \cdot \mb I({\mathcal{E}}_j)}
       \right)\mathbb{I}(\cup^N_{i=1}{\mathcal{E}}_i) 
       \Bigg\vert \nonumber \\
       &+
       \mathbb{I}(\cap^N_{i=1}{\mathcal{E}}^{^c}_i)  \label{eq:eq_A4} 
       \\
        =& \Bigg\vert
       \sum^N_{i=1} \left( 
       {\EE}_{\bm \gamma} \left[ Q_{t+1}(\bm x_{t}, \widetilde{\bm \xi}_{t+1})  \right] - 
       {\EE}_{\widetilde{ \bm \gamma}_{i}} \left[ Q_{t+1}(\bm x_{t}, \widetilde{\bm \xi}_{t+1}) \right]
       \right)\
       \frac{k(\Vert
       \bm \gamma-\widetilde{\bm \gamma}_{i}
       \Vert) \cdot \mb I({\mathcal{E}}_i)}{\sum^N_{j=1}k(\Vert\bm \gamma-\widetilde{\bm \gamma}_{j}
       \Vert) \cdot \mb I({\mathcal{E}}_j)}
       \mathbb{I}(\cup^N_{i=1}{\mathcal{E}}_i) \Bigg\vert  \nonumber
       \\
       & +
       \mathbb{I}(\cap^N_{i=1}{\mathcal{E}}^{^c}_i)  \nonumber
       \\
       \leq &
       \sum^N_{i=1} \left( 
       M_{t+1}\Vert \bm\gamma -\widetilde{ \bm \gamma}_{i}\Vert
       \cdot
        \Bigg\vert
       \frac{k(\Vert
       \bm \gamma-\widetilde{\bm \gamma}_{i}
       \Vert) \cdot \mb I({\mathcal{E}}_i)}{\sum^N_{j=1}k(\Vert\bm \gamma-\widetilde{\bm \gamma}_{j}
       \Vert) \cdot \mb I({\mathcal{E}}_j)}
        \mathbb{I}(\cup^N_{i=1}{\mathcal{E}}_i) \Bigg\vert\right) 
       +
       \mathbb{I}(\cap^N_{i=1}\mathcal{E}^{^c}_i)  \nonumber
       \\
       \leq & M_{t+1} h + \mathbb{I}(\cap^N_{i=1}\mathcal{E}^{^c}_i).
    \end{align}
    In \eqref{eq:def_m}, we expand the definition of \(\widetilde{m}(\cdot)\) over the complementary events \(\cup_{i=1}^N {\mathcal{E}}_i\) and \(\cap_{i=1}^N {\mathcal{E}}_i^{^{c}}\).
   The inequality \eqref{eq:eq_A4} holds since ${\mathbb{E}_{\bm\gamma}}[Q_{t+1}(\bm{x}_{t}, \widetilde{\bm{\xi}}_{t+1})] \leq 1$ by \ref{assumption:lipschitz_value}. Note that the indicator function $\mathbb{I}(\cap^N_{i=1}{\mathcal{E}}^{^c}_i)$ in the last inequality is a Bernoulli random variable with the parameter being the probability $\mathbb{P}(\cap^N_{i=1}{\mathcal{E}}^{^c}_i)$. Hence, we have $(a)\leq M_{t+1}h$ with probability $1-\mathbb{P}(\cap^N_{i=1}{\mathcal{E}}^{^c}_i)$. 
    Furthermore, we can establish an upper bound on the probability $\mathbb{P}(\cap^N_{i=1}{\mathcal{E}}^{^c}_i)$ by recursion. To maintain clarity and avoid overly cumbersome notation, assume that $N=4$ without loss of generality. Then, we have 
      \begin{align}
        \mathbb{P}(\cap^4_{i=1}{\mathcal{E}}^{^c}_i)
        =&\int_{\Xi_h^{^c}(\bm\gamma)} 
        \int_{\Xi_h^{^c}(\bm\gamma)} 
        \int_{\Xi_h^{^c} (\bm\gamma) } 
        \int_{\Xi_h^{^c} (\bm\gamma)}
         \mu (\mathrm{d}\bm \gamma_1)
        f (\bm \gamma_2 | \bm \gamma_1  ) 
        f (\bm \gamma_3 | \bm \gamma_2  ) 
        f (\bm \gamma_4 | \bm \gamma_3  ) 
        \mathrm{d} \bm \gamma_1 \mathrm{d} \bm \gamma_2 \mathrm{d} \bm \gamma_3 \mathrm{d} \bm \gamma_{4}  \nonumber
        \\
        =&
        \int_{\Xi_h^{^c}(\bm\gamma)} 
        \int_{\Xi_h^{^c}(\bm\gamma)} 
        \int_{\Xi_h^{^c} (\bm\gamma) } 
         \mu (\mathrm{d}\bm \gamma_1)
        f (\bm \gamma_2 | \bm \gamma_1  ) 
        f (\bm \gamma_3 | \bm \gamma_2  ) 
         \Big( 
        1 -\int_{\Xi_h (\bm\gamma)}
        f (\bm \gamma_{4} | \bm \gamma_{3}  ) 
        \Big)
        \mathrm{d} \bm \gamma_1 \mathrm{d} \bm \gamma_2 \mathrm{d} \bm \gamma_3 \mathrm{d} \bm \gamma_{4} 
        \label{eq:10}
        \\
        \leq&
        \int_{\Xi_h^{^c}(\bm\gamma)} 
        \int_{\Xi_h^{^c}(\bm\gamma)} 
        \int_{\Xi_h^{^c} (\bm\gamma) } 
         \mu (\mathrm{d}\bm \gamma_1)
        f (\bm \gamma_2 | \bm \gamma_1  ) 
        f (\bm \gamma_3 | \bm \gamma_2  ) 
        \mathrm{d} \bm \gamma_1 \mathrm{d} \bm \gamma_2 \mathrm{d} \bm \gamma_3 
        \cdot\Big(1-C_p h^p \underline{f}\Big)
        \label{eq:11}
        \\
        =&
        \int_{\Xi_h^{^c}(\bm\gamma)} 
        \int_{\Xi_h^{^c}(\bm\gamma)} 
         \mu (\mathrm{d}\bm \gamma_1)
        f (\bm \gamma_2 | \bm \gamma_1  ) 
        \Big( 
        1 -
        \int_{\Xi_h (\bm\gamma)} f (\bm \gamma_{3} | \bm \gamma_{2}  ) 
        \Big)
        \mathrm{d} \bm \gamma_1 \mathrm{d} \bm \gamma_2 \mathrm{d} \bm \gamma_3 
        \cdot\Big(1- C_p h^p \underline{f}\Big)
        \nonumber
                \\
        \leq&
        \int_{\Xi_h^{^c}(\bm\gamma)} 
        \int_{\Xi_h^{^c}(\bm\gamma)} 
         \mu (\mathrm{d}\bm \gamma_1)
        f (\bm \gamma_2 | \bm \gamma_1  ) 
        \mathrm{d} \bm \gamma_1 \mathrm{d} \bm \gamma_2 
        \cdot\Big(1 - C_p h^p\underline{f}\Big)^2
        \nonumber
        \\
        &\qquad\vdots \nonumber
        \\
        \leq&\left(1- C_ph^p \underline{f}\right)^{4},
        \label{eq:13}
    \end{align}
Here, for simplicity, we define $C_p = \C$, a constant that depends only on $p$. Specifically, $C_p$ arises in the closed-form expression for the volume of a $p$-dimensional Euclidean ball~(\citet[Equation 5.19.4]{lozier2003nist}). The function $\Gamma(\cdot)$, known as the gamma function, generalizes the factorial function. For positive integers $z$, $\Gamma(z) = (z-1)!$, while for non-integer inputs, it is an extension of the factorial.
In the equality \eqref{eq:10}, instead of directly integrating $f(\bm{\gamma}_4 | \bm{\gamma}_3)$ over ${\Xi_h^{^c}(\bm\gamma)}$, we integrate $f(\bm{\gamma}_4 | \bm{\gamma}_3)$ over $\Xi_h(\bm{\gamma}_t)$ and take the complimentary probability. 
Then, as in the inequality \eqref{eq:11}, we can extract the upper bound on the last integration by utilizing the lower bound $\underline{f}$ on the density function in \ref{a2}. 
By applying recursion, we arrive at \eqref{eq:13}.
Therefore, for any $N\in\mathbb{Z}_{++}$ we have 
\begin{equation}
         \mathbb{P}\left(
         (a) \leq M_{t+1}h 
         \right) \geq 1-(1-C_ph^p \underline{f})^{N}.
         \label{eq:bound_a}
\end{equation}

\noindent\underline{$(b)$:} 
Using the definitions in \eqref{eq:naive_E_hat} and \eqref{eq:function_m}, the term $(b)$ in \eqref{eq:decomposition} can be rewritten as follows:
 \begin{align}
       (b)=& 
         \Bigg\vert 
         \left(
         \sum_{i=1}^N 
         \left(
         \EE_{\widetilde{ \bm \gamma}_{i}}\left[ Q_{t+1}(\bm x_{t}, \widetilde{\bm \xi}_{t+1})  \right] 
         - 
         Q_{t+1}(\bm x_{t}, \widetilde{ \bm \gamma}_{i+1}) 
         \right)
        \frac{k(\Vert
       \bm \gamma-\widetilde{\bm \gamma}_{i}
       \Vert) \cdot \mb I({\mathcal{E}}_i)}{\sum^N_{j=1}k(\Vert\bm \gamma-\widetilde{\bm \gamma}_{j}
       \Vert) \cdot \mb I({\mathcal{E}}_j)}
         \right)\cdot\mathbb{I}(\cup^N_{i=1}{\mathcal{E}}_i) \nonumber
           \\
         &+(0-0)\cdot\mathbb{I}(\cap^N_{i=1}{\mathcal{E}}_i^{^c})\Bigg\vert 
         \nonumber
         \\
         =&\Bigg\vert 
         \left(
         \sum_{i=1}^N 
         \left(
         \EE_{\widetilde{ \bm \gamma}_{i}}\left[ Q_{t+1}(\bm x_{t}, \widetilde{\bm \xi}_{t+1})\right] 
         - 
         Q_{t+1}(\bm x_{t}, \widetilde{ \bm \gamma}_{i+1}) 
         \right)
                \frac{k(\Vert
       \bm \gamma-\widetilde{\bm \gamma}_{i}
       \Vert) \cdot \mb I({\mathcal{E}}_i)}{\sum^N_{j=1}k(\Vert\bm \gamma-\widetilde{\bm \gamma}_{j}
       \Vert) \cdot \mb I({\mathcal{E}}_j)}
         \right)\cdot\mathbb{I}(\cup^N_{i=1}{\mathcal{E}}_i)\Bigg\vert
         \nonumber
         \\
         =&
         \Bigg\vert 
         \sum_{i=1}^N 
         \scalebox{1.5}{$\Bigg($}
         \frac{\EE_{\widetilde{ \bm \gamma}_{i}}\left[ Q_{t+1}\left(\bm x_{t}, \phi( 
         \widetilde{\bm \gamma}_{i},\widetilde{\zeta}_{i+1}
         ) \right) \right] 
         - 
    Q_{t+1}
            \left(
                \bm x_{t}, \phi( 
                 \widetilde{\bm \gamma}_{i},\widetilde{\zeta}_{i+1}
                 ) 
             \right)}{k(\Vert
               \bm \gamma-\widetilde{\bm \gamma}_{i}
               \Vert)^{-1}}
         \scalebox{1.5}{$\Bigg)$}
         \Bigg\vert
         \frac{\mathbb{I}(\cup^N_{i=1}{\mathcal{E}}_i)}{\sum^N_{j=1}k(\Vert\bm \gamma-\widetilde{\bm \gamma}_{j}
       \Vert) \cdot \mb I({\mathcal{E}}_j)}.
         \label{eq:pre_bernstein}
    \end{align}
In \eqref{eq:pre_bernstein}, we replace the random variables $\widetilde{\bm\xi}_{t+1}$ and $\widetilde{\bm\gamma}_{i+1}$ with $\phi(\widetilde{\bm \gamma}_{i},\widetilde{\zeta}_{i+1})$ by Definition \ref{def:markov}. Note that, by conditioning on $\widetilde{\bm\gamma}_i={\bm\gamma}_i$ for all \(i \in [N]\), we can treat the terms $Q_{t+1}(\cdot,\phi({\bm \gamma}_{i},\widetilde{\zeta}_{2}))$, $\ldots$, $Q_{t+1}(\cdot,\phi({\bm \gamma}_{i},\widetilde{\zeta}_{N+1}))$ in \eqref{eq:pre_bernstein} as \emph{independent} random variables, since $\widetilde{\zeta}_2,\ldots,\widetilde{\zeta}_{N+1}$ are independent.
 From this, it is clear that \eqref{eq:naive_E_hat} is a sum of independent random variables, while \eqref{eq:function_m} represents the expectation of this sum.
 Hence, conditional on $\widetilde{\bm{\gamma}}_{[1:N]}=\bm{\gamma}_{[1:N]}$ and the event $\cup^N_{i=1}{\mathcal{E}}_i$, we can apply Bernstein inequality~(\citet[Theorem 2.8.4]{vershynin2018high}) to \eqref{eq:pre_bernstein}, as follows: 
 \begin{align}
    &\mathbb{P}\left(
    (b) \geq \epsilon \; \middle\vert \; \widetilde{\bm{\gamma}}_{[1:N]}=\bm{\gamma}_{[1:N]}, \;\cup^N_{i=1}{\mathcal{E}}_i
    \right)
    \\
    =& 2\cdot
    \mathbb{P}\left(
    \sum_{i:\Vert\bm \gamma - \bm\gamma_{i}\Vert\leq h}
     \scalebox{1.5}{$\Bigg($}
        \frac{\EE_{{ \bm \gamma}_{i}}\left[ Q_{t+1}\left(\bm x_{t}, \phi( 
         {\bm \gamma}_{i},\widetilde{\zeta}_{i+1}
         ) \right)  \right] 
         - 
         Q_{t+1}\left(\bm x_{t}, \phi( 
         {\bm \gamma}_{i},\widetilde{\zeta}_{i+1}
         ) \right)}{k(\Vert\bm \gamma-{\bm \gamma}_{i}\Vert)^{-1}}
    \scalebox{1.5}{$\Bigg)$}
         \geq 
         \epsilon 
         \sum^N_{j=1}k(\Vert\bm \gamma-{\bm \gamma}_{j}
       \Vert) \cdot \mb I({\mathcal{E}}_j)
    \right) 
    \label{eq:bernstein_placeholder}\\
    \leq & 2\cdot \exp\left(
    -\frac{1}{2} \cdot
    \frac{
    \left(
        \epsilon 
         \sum^N_{j=1}k(\Vert\bm \gamma-{\bm \gamma}_{j}
       \Vert) \cdot \mb I({\mathcal{E}}_j)
    \right)^2}
    {\sum^N_{i=1}\VV_{\bm \gamma_{i}} \left[Q_{t+1}\left(\bm x_{t}, \phi( 
         {\bm \gamma}_{i},\widetilde{\zeta}_{i+1}
         )\right) \right] \mathbb{I}({\mathcal{E}}_i)
         \cdot
         k(\Vert\bm \gamma-{\bm \gamma}_{i}
       \Vert)^2 + \frac{1}{3}
       \left(
        \epsilon 
         \sum^N_{j=1}k(\Vert\bm \gamma-{\bm \gamma}_{j}
       \Vert) \cdot \mb I({\mathcal{E}}_j)
    \right)
       }
    \right) \\
    =& 2\cdot \exp\left(
    -\frac{\epsilon^2}{2} \cdot
    \frac{
        \sum^N_{j=1}k(\Vert\bm \gamma-{\bm \gamma}_{j}
               \Vert) \cdot \mb I({\mathcal{E}}_j)
    }
    {\sum^N_{i=1}\VV_{\bm \gamma_{i}} \left[ Q_{t+1}\left(\bm x_{t}, \phi( 
         {\bm \gamma}_{i},\widetilde{\zeta}_{i+1}
         )\right) \right]\frac{
         k(\Vert\bm \gamma-{\bm \gamma}_{i}
       \Vert)^2\cdot \mathbb{I}({\mathcal{E}}_i)}{\sum^N_{j=1}k(\Vert\bm \gamma-{\bm \gamma}_{j}
       \Vert)\cdot \mb I({\mathcal{E}}_j)} 
         + \frac{\epsilon}{3} }
    \right).
\end{align}
 By the tower rule, taking expectations on both sides yields 
\begin{align}
    \mathbb{P}\left(
    (b) \geq \epsilon
    \right)
    &
    \leq 
   2\cdot \EE \left[
     \exp\left(
    -\frac{\epsilon^2}{2} \cdot
    \frac{
        \sum^N_{j=1}k(\Vert\bm \gamma-\widetilde{\bm \gamma}_{j}
               \Vert) \cdot \mb I({\mathcal{E}}_j)
    }
    {\sum^N_{i=1}\VV_{\widetilde{\bm \gamma}_{i}} \left[ Q_{t+1}\left(\bm x_{t}, \phi( 
         {\bm \gamma}_{i},\widetilde{\zeta}_{i+1}
         )\right)\right]\frac{
         k(\Vert\bm \gamma-\widetilde{\bm \gamma}_{i}
       \Vert)^2\cdot \mathbb{I}({\mathcal{E}}_i)}{\sum^N_{j=1}k(\Vert\bm \gamma-\widetilde{\bm \gamma}_{j}
       \Vert)\cdot \mb I({\mathcal{E}}_j)} 
         + \frac{\epsilon}{3} }
    \right)
    {\mathbb{I}(\cup^N_{i=1}{\mathcal{E}}_i)}
    \right] 
    \nonumber \\
    &\leq 2\cdot \EE \left[
    \exp\left(
    -\frac{\epsilon^2}{2} \cdot
    \frac{
    \underline{k}
         {\sum^N_{j=1}\mathbb{I}({\mathcal{E}}_j)}
    }
    {3\overline{k}M_{t+1}h+\overline{k}\VV \left[ Q_{t+1}(\bm x_{t}, \widetilde{\bm \xi}_{t+1}) \middle\vert \bm \gamma \right]  + \frac{\epsilon}{3} }
    \right)
    {\mathbb{I}(\cup^N_{i=1}{\mathcal{E}}_i)}
    \right]    \label{eq:first_ineq_variance}
    \\
    &= 2\cdot \EE \left[
    \exp\left(   -\epsilon^2 \beta(\epsilon)
    {\sum^N_{j=1}\mathbb{I}({\mathcal{E}}_j)}
    \right)
    {\mathbb{I}(\cup^N_{i=1}{\mathcal{E}}_i)}
    \right] \label{eq:bound_simple} \\
    &
    \leq
    2\cdot
    \left(
    \EE \left[
    \exp\left(   -\epsilon^2 \beta(\epsilon)
    {\sum^N_{j=1}\mathbb{I}({\mathcal{E}}_j)}
    \right)
    \right]
    - 
    \mathbb{P}(\cap^N_{i=1}{\mathcal{E}}_i^{^c})
    \right) \label{eq:mean_exp}
    \\
    &\leq 
    2\cdot \EE \left[
    \exp\left(   -\epsilon^2 \beta(\epsilon)
    {\sum^N_{j=1}\mathbb{I}({\mathcal{E}}_j)}
    \right)
    \right] 
    -
    2\cdot(1-C_p h^p \underline{f})^{N}
    , 
    \label{eq:bound_b}
\end{align}
 where the constant $\underline{k} \geq 0$ depends only on the choice of the function $k(\cdot)$ in \eqref{eq:kernel}, and, in \eqref{eq:bound_simple}, we simplify the terms inside the exponential function by defining $\beta(\epsilon) = \frac{1}{2} \cdot \frac{\underline{k}}{3 \overline{k} M_{t+1} h + \overline{k} \VV_{\bm{\gamma}}[ Q_{t+1}(\bm{x}_t, \widetilde{\bm{\xi}}_{t+1})] + \frac{\epsilon}{3}}$.
Note that the indicator functions \(\mathbb{I}({\mathcal{E}}_j)\) for all \(j \in [N]\), and \(\mathbb{I}(\cup_{i=1}^N {\mathcal{E}}_i)\) are random variables as we take expectations. 
In \eqref{eq:first_ineq_variance}, we bound the denominator inside the exponential function by applying both Lemma \ref{lem:lip_variance} and \ref{assumption_kernel}: conditioning on $\cup^N_{i=1}{\mathcal{E}}_i$, for any $\bm \gamma_{[1:N]}$ and $\bm \gamma$, we have
\begin{align}
    & {\sum^N_{i=1}
    \VV_{{\bm \gamma}_{i}} \left[ Q_{t+1}\left(\bm x_{t}, \phi( 
         {\bm \gamma}_{i},\widetilde{\zeta}_{i+1}
         )\right)   \right] 
         \frac{
         k(\Vert\bm \gamma-{\bm \gamma}_{i}
       \Vert)^2\cdot \mathbb{I}({\mathcal{E}}_i)}{\sum^N_{j=1}k(\Vert\bm \gamma-{\bm \gamma}_{j}
       \Vert)\cdot \mb I({\mathcal{E}}_j)} 
         }
     \nonumber
    \\
    =  &{\sum^N_{i=1}\VV_{{\bm \gamma}_{i}} \left[ Q_{t+1}(\bm x_{t}, \widetilde{\bm \gamma}_{i+1})\right] 
    \frac{
         k(\Vert\bm \gamma-{\bm \gamma}_{i}
       \Vert)^2\cdot \mathbb{I}({\mathcal{E}}_i)}{\sum^N_{j=1}k(\Vert\bm \gamma-{\bm \gamma}_{j}
       \Vert)\cdot \mb I({\mathcal{E}}_j)} 
    }
    \nonumber
    \\
    \leq& \overline{k} \left(
    {\sum^N_{i=1} 
 \VV_{{\bm \gamma}_{i}} \left[ Q_{t+1}(\bm x_{t}, \widetilde{\bm \gamma}_{i+1})\right] 
    \frac{
         k(\Vert\bm \gamma-{\bm \gamma}_{i}
       \Vert)\cdot \mathbb{I}({\mathcal{E}}_i)}{\sum^N_{j=1}k(\Vert\bm \gamma-{\bm \gamma}_{j}
       \Vert)\cdot \mb I({\mathcal{E}}_j)} 
    - \VV_{\bm \gamma} \left[ Q_{t+1}(\bm x_{t}, \widetilde{\bm \xi}_{t+1})  \right]
    +
    \VV_{\bm \gamma} \left[ Q_{t+1}(\bm x_{t}, \widetilde{\bm \xi}_{t+1}) \right]
    } 
       \right)
   \nonumber
 \\
    \leq& \left\vert \overline{k} 
    \sum^N_{i=1} \left(
 \VV_{{\bm \gamma}_{i}} \left[ Q_{t+1}(\bm x_{t}, \widetilde{\bm \gamma}_{i+1}) \right] 
  - \VV_{\bm \gamma} \left[ Q_{t+1}(\bm x_{t}, \widetilde{\bm \xi}_{t+1})  \right]
    \right)\cdot
        \frac{
         k(\Vert\bm \gamma-{\bm \gamma}_{i}
       \Vert)\cdot \mathbb{I}({\mathcal{E}}_i)}{\sum^N_{j=1}k(\Vert\bm \gamma-{\bm \gamma}_{j}
       \Vert)\cdot \mb I({\mathcal{E}}_j)} 
    \right\vert
    \nonumber
    \\
    &+
    \left\vert
    \overline{k}
    \VV_{\bm \gamma} \left[ Q_{t+1}(\bm x_{t}, \widetilde{\bm \xi}_{t+1}) \right]
    \right\vert   
    \nonumber
       \\
    \leq  & 3\overline{k}M_{t+1}h+\overline{k}\VV_{\bm \gamma} \left[ Q_{t+1}(\bm x_{t}, \widetilde{\bm \xi}_{t+1})\right].  
    \label{eq:lip_var}
\end{align}
Here, the inequality~\eqref{eq:lip_var} holds due to the Lipschitz continuity of the conditional variance established in \Cref{lem:lip_variance}.
  In \eqref{eq:mean_exp}, we apply the decomposition of the following expectation using the tower rule: 
  \begin{align*}
   \EE \left[
    \exp\left(   -\epsilon^2 \beta(\epsilon)
    {\sum^N_{j=1}\mathbb{I}({\mathcal{E}}_j)}
    \right)
    \right]
    =&
    \EE \left[
    \exp\left(   -\epsilon^2 \beta(\epsilon)
    {\sum^N_{j=1}\mathbb{I}({\mathcal{E}}_j)}
    \right) \middle\vert
    \mathbb{I}(\cup^N_{i=1}{\mathcal{E}}_i)=1
    \right]\cdot
    \mathbb{P}( \cup^N_{i=1}{\mathcal{E}}_i)
    \\ &
    +
    \exp\left(0
    \right)  
    \cdot
    \mathbb{P}( \cap^N_{i=1}{\mathcal{E}}^{^{c}}_i).
  \end{align*}
 Note that \(\sum_{j=1}^N \mathbb{I}({\mathcal{E}}_j)\) can be considered as a sum of \(N\) \emph{dependent} Bernoulli random variables. Given this, the first term in \eqref{eq:bound_b}, i.e., the expectation, can be explicitly written, and its upper bound can be established:
\begin{align}
    & \EE \left[
    \exp\left(   -\epsilon^2 \beta(\epsilon)
    {\sum^N_{j=1}\mathbb{I}({\mathcal{E}}_j)}
    \right)
    \right] \nonumber \\
    =&\EE \left[
    \exp\left(   -\epsilon^2 \beta(\epsilon)
    \mathbb{I}({\mathcal{E}}_1)
    \right)\cdots \exp\left(   -\epsilon^2 \beta(\epsilon)
    \mathbb{I}({\mathcal{E}}_N)
    \right)
    \right] \nonumber \\
   =&
 \int_{\Xi} \int_{\Xi}  \cdots \int_{\Xi} \int_{\Xi}\left(\exp \left(-\epsilon^2 \beta(\epsilon) \mathbb{I}\left({\mathcal{E}}_1\right)\right) \mu\left(\mathrm{d} \bm\gamma_1\right)\right)   \nonumber \\
&\qquad\qquad\qquad\;\;
\left(\exp \left(-\epsilon^2 \beta(\epsilon) \mathbb{I}\left({\mathcal{E}}_2\right)\right) f\left(\bm\gamma_2 \mid \bm\gamma_1\right)\right)\cdots \nonumber \\
&\qquad\qquad\qquad\;\;
\left(\exp \left(-\epsilon^2 \beta(\epsilon) \mathbb{I}\left({\mathcal{E}}_N\right)\right) f\left(\bm\gamma_{N} \mid \bm\gamma_{N-1}\right)\right) \mathrm{d} \bm\gamma_1 \mathrm{d} \bm\gamma_2  \cdots \mathrm{d} \bm\gamma_{N-1} \mathrm{d} \bm\gamma_{N}. \nonumber \\
    \leq&\left(1+C_ph^p\underline{f}
    \left(
    \exp(-\epsilon^2\beta(\epsilon)) - 1
    \right)
    \right)^N. 
    \label{eq:bound_b_exp}
\end{align}
Here, we obtain \eqref{eq:bound_b_exp} by recursion similar to \eqref{eq:13}.
Combining \eqref{eq:bound_a} and \eqref{eq:bound_b}, with \eqref{eq:bound_b_exp} replacing the expectation in \eqref{eq:bound_b}, we have
\begin{align}
\mathbb{P}\left[ (a) +(b) \geq M_{t+1}h+\epsilon \right] 
&\leq 
2\cdot \left(1+C_ph^p\underline{f}
    \left(
    \exp(-\epsilon^2\beta(\epsilon)) - 1
    \right)
    \right)^N - 2\cdot(1-C_ph^p \underline{f})^{N} + (1-C_ph^p \underline{f})^{N} \nonumber
    \\
    & \leq 2\cdot \left(1+C_ph^p\underline{f}
    \left(
    \exp(-\epsilon^2\beta(\epsilon)) - 1
    \right)
    \right)^N \nonumber
    \\
    &\leq 2\cdot \exp\left(
    NC_ph^p\underline{f}
    \left(
    \exp(-\epsilon^2\beta(\epsilon)) - 1
    \right)
    \right) 
    \label{eq:(1+x)_ineq}
    \\
    & \leq
    2\cdot \exp\left(
    NC_ph^p\underline{f}
    \left( \frac{-\epsilon^2\beta(\epsilon)}{\alpha}
    \right)
    \right)  \label{eq:custom_ineq}
    \\
    & =
    2\cdot \exp\left(
    -\frac{N C_p h^p \underline{f}\underline{k}}{2\alpha} \cdot
    \frac{\epsilon^2}{\overline{k}\VV \left[ Q_{t+1}(\bm x_{t}, \widetilde{\bm \xi}_{t+1}) \middle\vert \bm \gamma \right] + 3\overline{k}M_{t+1}h  + \frac{1}{3}\epsilon }
    \right), \label{eq:finite_bound}
\end{align}
where $\alpha>0$ is a sufficiently large constant. 
The inequality \eqref{eq:(1+x)_ineq} holds since $(1+x)^N\leq \exp(Nx)$ holds for any $x\geq-1$ and any $N\in\mathbb{Z}_{+}$. 
In \eqref{eq:custom_ineq}, we use the following fact: under the boundedness assumption on $\Xi$ in \ref{a2}, for any finite $\epsilon\in[0,1]$, there exists a universal constant $\alpha$ such that 
\begin{equation*}
    \exp\left(-\epsilon^2\beta(\epsilon)\right) -1 \leq - \frac{\epsilon^2\beta(\epsilon)}{\alpha}.
\end{equation*}
This inequality can be verified by checking the first derivative condition of $\exp(-\epsilon^2\beta(\epsilon))$. This completes the proof.
\end{proof}

\noindent \Cref{prop:kernel_estimation_error} asserts that as the sample size $N$ increases, the probability that the estimator deviates from the true conditional expectation decreases exponentially fast. By setting \eqref{eq:prob_kernel_estimation_errors} equal to $\delta$, we derive the following inequality for a \emph{fixed} decision $\bm{x}_t$.
\begin{corollary}
\label{corol:bound_for_a_fixed}
    For any fixed $\delta\in(0,1]$, $h\geq0$, and $\bm x_t\in\mc X_t(\bm x_{t-1},\bm \xi_{t})$, we have   
    \begin{multline}
        \left| \EE_{\bm \gamma} \left[ Q_{t+1}(\bm x_{t}, \widetilde{\bm \xi}_{t+1}) \right] - \widehat \EE_{\bm \gamma} \left[ Q_{t+1}(\bm x_{t}, \widetilde{\bm \xi}_{t+1})\right] \right| \\ 
        \leq M_{t+1} h + \frac{2\alpha\log\left(\frac{2}{\delta}\right)}{3NC_ph^p\underline{f}\underline{k}}+\sqrt{ \frac{2\alpha\overline{k}\log\left(\frac{2}{\delta}\right)\left(
        \sigma^2_{t+1}
        +3M_{t+1}h
        \right)}{NC_ph^p\underline{f}\underline{k}}}
    \end{multline}
   with probability at least $1-\delta$. Here, all other parameters are defined as in \Cref{prop:kernel_estimation_error}.
\end{corollary}
\begin{proof}
Setting $\delta$ equal to  \eqref{eq:finite_bound}, we have
\begin{equation}
    3\left(\frac{NC_ph^p\underline{f}\underline{k}}{2\alpha}\right)\epsilon^2-\log\left( \frac{2}{\delta} \right) \epsilon
    -3\log\left( \frac{2}{\delta} \right)\left(
    \overline{k}\VV_{\bm \gamma} \left[ Q_{t+1}(\bm x_{t}, \widetilde{\bm \xi}_{t+1})  \right] +3\overline{k}M_{t+1} h
    \right)=0.
    \label{eq:quad_eq}
\end{equation}
Then, we can solve the quadratic equation \eqref{eq:quad_eq} for $\epsilon$ and further establish an upper bound on $\epsilon$ as follows:
\begin{align}
    \epsilon &= \frac{\log\left(\frac{2}{\delta} \right)+\sqrt{\left(\log\left( \frac{2}{\delta} \right)\right)^2+36\left(\frac{NC_ph^p\underline{f}\underline{k}}{2\alpha}\right)\log\left( \frac{2}{\delta} \right)\left(
    \overline{k}\VV_{\bm \gamma} \left[ Q_{t+1}(\bm x_{t}, \widetilde{\bm \xi}_{t+1})  \right] +3\overline{k}M_{t+1} h
    \right)}}{6\left(\frac{NC_ph^p\underline{f}\underline{k}}{2\alpha}\right)}
    \nonumber
    \\
    &\leq 
    \frac{2\log\left(\frac{2}{\delta} \right)+\sqrt{36\left(\frac{NC_ph^p\underline{f}\underline{k}}{2\alpha}\right)\log\left( \frac{2}{\delta} \right)\left(\overline{k}
    \VV_{\bm \gamma} \left[ Q_{t+1}(\bm x_{t}, \widetilde{\bm \xi}_{t+1})  \right] +3\overline{k}M_{t+1} h
    \right)}}{6\left(\frac{NC_ph^p\underline{f}\underline{k}}{2\alpha}\right)}
    \nonumber
    \\
    &= \frac{2\alpha\log\left(\frac{2}{\delta}\right)}{3NC_ph^{^p}\underline{f}\underline{k}}+\sqrt{ \frac{2\alpha\overline{k}\log\left(\frac{2}{\delta}\right)\left(
        \VV_{\bm \gamma}\left[ Q_{t+1}(\bm x_{t}, \widetilde{\bm \xi}_{t+1})  \right]+3M_{t+1}h
        \right)}{NC_ph^{^p}\underline{f}\underline{k}} }.
    \label{eq:simplified_bound}
\end{align}
The inequality above holds since $ x+\sqrt{x^2+y}\leq 2x+\sqrt{y}$ for any $x,y\in\mathbb{R}_{+}$. By simplifying \eqref{eq:simplified_bound} using the upper bound $\sigma^2_{t+1}$ on the conditional variance from \Cref{a7}, we obtain the claim. 
\end{proof}
\noindent We now extend this result to establish the uniform generalization bound for all decisions $\bm{x}_t$ within the feasible region.

\begin{proposition}\label{thm:Generalization Bound_for_a_Continuous}
\label{lem:generalization_bound_continuousX} 
Suppose that $D_t$ is a finite diameter from \ref{a3}, and $L_{t+1}$ is the Lipschitz constant from \ref{assumption:lipschitz_value}. Then, for any fixed $\delta\in(0,1]$, $h\geq0$, $\eta > 0$ and $\bm \gamma \in \Xi$, we have
\begin{align*}
    &\left| \EE_{\bm \gamma} \left[ Q_{t+1}(\bm x_{t}, \widetilde{\bm \xi}_{t+1})  \right] - \widehat \EE_{\bm \gamma} \left[ Q_{t+1}(\bm x_{t}, \widetilde{\bm \xi}_{t+1}) \right] \right| 
    \\ &\leq 
        M_{t+1} h + \frac{2\alpha
        \log \Big(\frac{
         \mathcal{O}(1) \left({D_{t}}/{\eta}\right)^{d_{t}} 
        }{\delta}\Big)
        }{3NC_ph^p\underline{f}\underline{k}}+\sqrt{ \frac{
        2\alpha\overline{k}
        \log \Big(\frac{
         \mathcal{O}(1) \left({D_{t}}/{\eta}\right)^{d_{t}} 
        }{\delta}\Big)
        \left(
        \sigma^2_{t+1}+3M_{t+1}h
        \right)}{NC_ph^p\underline{f}\underline{k}} } + 2 L_{t+1} \eta,  
        \\
        & \qquad\qquad\qquad\qquad\qquad\qquad
        \qquad\qquad\qquad\qquad\qquad\qquad
        \qquad\qquad\qquad\;\;\;
        \forall \x_{t}\in \mc X_{t}(\bm x_{t-1}, \bm \gamma )
\end{align*}
with probability at least $1-\delta$. Here, all other parameters are defined as in \Cref{prop:kernel_estimation_error}.
\end{proposition}

\begin{proof}
\noindent We define an $\eta$-net of the feasible region, i.e., a set of finite points $\mc X_{t}^{\eta}\left(\bm x_{t-1}, \bm\gamma \right) \subset \mc X_{t}\left(\bm x_{t-1}, \bm \gamma\right)$ such that for any $\x_{t} \in \mc X_{t}\left(\bm x_{t-1}, \bm \gamma \right)$, there exists $\bm y_{t}\in \mc X_{t}^{\eta}\left(\bm x_{t-1}, \bm\gamma \right)$ such that $\Vert \x_{t} - \bm y_{t} \Vert \leq \eta$. According to \citet{shapiro2005complexity}, its cardinality is $\vert \mc X_{t}^{\eta}\left(\bm x_{t-1}, \bm \gamma \right) \vert=\mathcal{O}(1) (D_{t}/\eta)^{d_{t}}$, where $D_t$ is the finite diameter as stated in \ref{a3}, $d_t$ is the dimension of $\mathcal{X}_t(\cdot)$, and $\eta > 0$ is a parameter that controls the number of discretized points. Since the value function $Q_{t+1}(\bm x_{t}, \bm \xi_{t+1})$ is $L_{t+1}$-Lipschitz continuous in $\x_{t}$, \Cref{corol:expect_lipschitz} implies that for any $\x_{t} \in \mc X_{t}\left(\bm x_{t-1}, \bm \gamma \right)$, there exists $\bm y_{t}\in \mc X_{t}^{\eta}\left(\bm x_{t-1}, \bm \gamma \right)$, $\Vert \x_{t} - \bm y_{t} \Vert \leq \eta$, such that:
\begin{equation}
\label{eq:lipschitz_lma_Gen}
\begin{aligned}
\left| \EE_{\bm \gamma} \left[ Q_{t+1}(\bm x_{t}, \widetilde{\bm \xi}_{t+1}) \right] - \EE_{\bm \gamma} \left[ Q_{t+1}(\bm y_{t}, \widetilde{\bm \xi}_{t+1}) \right] \right| & \leq L_{t+1} \eta.
\end{aligned}
\end{equation}
Furthermore, from \Cref{corol:bound_for_a_fixed}, for a fixed $\bm y_{t}\in \mc X_{t}^{\eta}(\bm x_{t-1}, \bm \gamma ),$ we have
\begin{multline}
        \left| \EE_{\bm \gamma} \left[ Q_{t+1}(\bm y_{t}, \widetilde{\bm \xi}_{t+1}) \right] - \widehat \EE_{\bm \gamma} \left[ Q_{t+1}(\bm y_{t}, \widetilde{\bm \xi}_{t+1}) \right] \right| \\ 
        \leq M_{t+1} h + \frac{2\alpha\log\left(\frac{2}{\delta}\right)}{3NC_ph^p\underline{f}\underline{k}}+\sqrt{ \frac{2\alpha\overline{k} \log \left(\frac{2}{\delta}\right)\left(
        \sigma^2_{t+1}+3M_{t+1}h
        \right)}{NC_ph^p\underline{f}\underline{k}} }
    \end{multline}
   with probability at least $1-\delta$.
Applying union bound over $\mc X_{t}^{\eta}(\cdot)$, we get
\begin{align*}
        &\left| \EE_{\bm \gamma} \left[ Q_{t+1}(\bm y_{t}, \widetilde{\bm \xi}_{t+1})\right] - \widehat \EE_{\bm \gamma} \left[ Q_{t+1}(\bm y_{t}, \widetilde{\bm \xi}_{t+1})\right] \right| \\ 
        \leq& M_{t+1} h + \frac{2\alpha
        \log \left(\frac{2\lvert \mc X_{t}^{\eta}\left(\bm x_{t-1}, \bm \gamma \right) \rvert}{\delta}\right)
        }{3NC_ph^p\underline{f}\underline{k}}+\sqrt{ \frac{2\alpha\overline{k}\log \left(\frac{2\lvert \mc X_{t}^{\eta}\left(\bm x_{t-1}, \bm \gamma \right) \rvert}{\delta}\right)\left(
        \sigma^2_{t+1}+3M_{t+1}h
        \right)}{NC_ph^p\underline{f}\underline{k}} } 
        \\  = & 
        M_{t+1} h + \frac{2\alpha
        \log \left(\frac{ \mc
         O(1) (D_{t}/\eta)^{d_{t}} 
        }{\delta}\right)
        }{3NC_ph^p\underline{f}\underline{k}}+\sqrt{ \frac{2\alpha\overline{k}
        \log \left(\frac{ \mc
         O(1) (D_{t}/\eta)^{d_{t}} 
        }{\delta}\right)
        \left(
        \sigma^2_{t+1}+3M_{t+1}h
        \right)}{NC_ph^p\underline{f}\underline{k}} },
        \\ &
        \qquad\qquad\qquad\qquad\qquad\qquad
        \qquad\qquad\qquad\qquad\qquad\qquad
        \qquad\qquad
        \forall \bm y_{t}\in \mc X_{t}^{\eta}(\bm x_{t-1}, \bm \gamma )
    \end{align*}
    with probability at least $1-\delta$.
Using the Lipschitz continuity of the expected value function in $\bm\x_{t}$ from \Cref{corol:expect_lipschitz}, we get
\begin{align*}
    &\left| \EE_{\bm \gamma} \left[ Q_{t+1}(\bm x_{t}, \widetilde{\bm \xi}_{t+1})\right] - \widehat \EE_{\bm \gamma} \left[ Q_{t+1}(\bm x_{t}, \widetilde{\bm \xi}_{t+1})\right] \right| 
    \\ \leq &
        M_{t+1} h + \frac{2\alpha
        \log \left(\frac{
         \mc O(1) (D_{t}/\eta)^{d_{t}} 
        }{\delta}\right)
        }{3NC_ph^p\underline{f}\underline{k}}+\sqrt{ \frac{2\alpha \overline{k}
        \log \left(\frac{
         \mc O(1) (D_{t}/\eta)^{d_{t}} 
        }{\delta}\right)
        \left(
        \sigma^2_{t+1}
        +
        3M_{t+1}h
        \right)}{NC_ph^p\underline{f}\underline{k}} } + 2 L_{t+1} \eta,  
        \\
        & \qquad\qquad\qquad\qquad\qquad
        \qquad\qquad\qquad\qquad\qquad\qquad
        \qquad\qquad\qquad\quad \;
        \forall \x_{t}\in \mc X_{t}(\bm x_{t-1}, \bm \gamma )
\end{align*}
with probability at least $1-\delta$. This completes the proof.
\end{proof}

\noindent So far, we have derived the generalization bound for an arbitrary stage $t$, focusing solely on the approximation error introduced by the estimator while treating the value function $Q_{t+1}(\cdot)$ as if it were available. As briefly discussed in \Cref{sec:problem_statement}, dynamic programming enables us to decompose the MSP problem into subproblems, where the value function for the MRST problem at each stage, $\widehat{Q}_{t+1}(\cdot)$, exhibits errors that accumulate from the approximation errors of the value functions at later stages, $[T] \setminus [t+1]$.
Therefore, by recursively applying \Cref{lem:generalization_bound_continuousX} from the terminal stage to the first, we derive the generalization bound for the MRST problem, as presented in the following theorem.

\begin{theorem}[Generalization Bound]\label{thm:out_of_sample}
\setlength{\abovedisplayskip}{0pt}%
Suppose that $\widetilde{\bm \gamma}^{(a)}_{[1:N+1]}$ and $\widetilde{\bm \gamma}^{(b)}_{[1:N+1]}$ represent two i.i.d.~Markovian trajectories, used in \eqref{eq:approx_dp_even} and \eqref{eq:approx_dp_odd}, respectively, to construct the MRST problem. Then, for any fixed $\delta_t \in (0,1]$ with $t\in[T]\setminus\{1\}$, $h\geq0$, $\eta>0$, and $\bm x_1 \in \mc X_1(\bm x_0, \bm \xi_1)$, we have
\begin{multline}
\label{eq_OOS}
        \bigg|    \EE_{\bm \xi_1} \big[ Q_2 (\bm x_1, \widetilde{\bm \xi}_{2} ) \big] 
        -
         \widehat \EE_{\bm \xi_1} \big[ \widehat Q_{2}(\bm x_1,  \xit_{2}) \big]  \bigg| \leq 
        \\
        \sum_{t=2}^{T}
        M_{t}h 
        +
        \frac{2\alpha\log \Big(\frac{\mc O(1) N^{t-2} \prod_{s =1}^{t-1} \left(\frac{D_{s}}{\eta}\right)^{d_s} }{ \delta_t }\Big)}{3NC_ph^p\underline{f}\underline{k}}
        +
         \sqrt{ \frac{2\alpha\overline{k}\log \Big(\frac{\mc O(1) N^{t-2} \prod_{s =1}^{t-1} \left(\frac{D_{s}}{\eta}\right)^{d_s} }{ \delta_t }\Big)\left(
        \sigma^2_{t}
        +3M_{t}h
        \right)}{NC_ph^p\underline{f}\underline{k}} }
        +
        2 L_t \eta
\end{multline}
with probability at least $1-\sum_{t=2}^{T} \delta_t$. Here, all other parameters are defined as in \Cref{prop:kernel_estimation_error} and \Cref{thm:Generalization Bound_for_a_Continuous}.
\end{theorem}
\begin{proof}
We begin with deriving the upper bound on the approximate value function errors at the terminal stage and accumulate the errors by backward induction. 
Without loss of generality, let us assume that the time horizon $T$ is an odd number---the only difference when $T$ is an even number is that the samples at the terminal stage are drawn from the second trajectory $\widetilde{\gamma}^{(b)}_{[1:N+1]}$, instead of the first trajectory $\widetilde{\gamma}^{(a)}_{[1:N+1]}$.

\noindent \underline{Stage $t = T$}: We have $Q_T(\bm x_{T-1}, \bm \xi_{T}) = \widehat Q_T(\bm x_{T-1}, \bm \xi_{T})$ for all combination of inputs $\bm x_{T-1}\in \mc X_{T-1}(\cdot)$ and $\bm \xi_T\in\Xi$ since $\EE_{\bm \xi_{T}}[Q_{T+1}(\bm x_{T},\widetilde{\bm \xi}_{T+1})]= 0$. To simplify the exposition, we omit the time indices on \(\bm x_t\) and \(\bm \xi_t\), using \(\bm x\) and \(\bm \xi\) instead. Additionally, we use \(\bm x'\) and \(\bm \xi'\) to distinguish \(\bm x\) and \(\bm \xi\) across different time stages.
From \Cref{lem:generalization_bound_continuousX}, for any fixed $\x'\in\mathcal{X}_{T-2}(\cdot)$ and $\bm\xi'\in\Xi$ we have
\begin{align*}
    &
    \left| \EE_{\bm \xi'} \left[ Q_T(\bm x,\xit) \right] - \widehat \EE_{\bm \xi'} \left[ \widehat Q_T(\bm x, \xit)\right] \right| 
    \\ &\leq 
        M_T h + \frac{2\alpha
        \log 
        \Bigg(
        \frac{
         \mc O (1) \left(\frac{D_{T-1}}{\eta}\right)^{d_{T-1}} 
        }{\delta_T}
        \Bigg)
        }{3NC_ph^p\underline{f}\underline{k}}+\sqrt{ \frac{2\alpha \overline{k}
        \log \Bigg(
        \frac{
         \mc O(1) \left(\frac{D_{T-1}}{\eta}\right)^{d_{T-1}} 
        }{\delta_T}
        \Bigg)
        \left(
        \sigma^2_T+3M_Th
        \right)}{NC_ph^p\underline{f}\underline{k}} } + 2 L_{T} \eta, 
         \\
        & \qquad\qquad\qquad\qquad\qquad\qquad
        \qquad\qquad\qquad\qquad\qquad\qquad
        \qquad\qquad\qquad\qquad\qquad\;
        \forall \bm x\in \mc X_{T-1}(\bm x', \bm \xi' )
\end{align*}
with probability at least $1-\delta_T$. Taking minimization over $\mc X_{T-1}(\bm x', \bm \xi' )$, for any fixed $\x'\in\mathcal{X}_{T-2}(\cdot)$ and $\bm\xi'\in\Xi$ we have
\begin{IEEEeqnarray*}{lCl}
  && \bigg| \min_{ \x \in \mc X_{T-1}\left(\bm x', \bm\xi' \right) } \left\lbrace 
  c_{_{T-1}}(\bm x,\bm\xi')
  + \EE_{\bm \xi'} \left[ Q_T(\bm x, \xit) \right]  \right\rbrace 
  - 
  \min_{ \x \in \mc X_{T-1}\left(\bm x', \bm\xi' \right) } \left\lbrace c_{_{T-1}}(\bm x,\bm\xi') + \widehat\EE_{\bm \xi'} \left[ \widehat Q_T(\bm x, \xit)\right]  \right\rbrace  
  \bigg| \\
 &=& \left|Q_{T-1}(\bm x',\bm \xi')  -\widehat Q_{T-1}(\bm x',\bm \xi') \right| \\
 & \leq &   M_T h + \frac{2\alpha
        \log \Bigg(\frac{
         \mc O(1) \left(\frac{D_{T-1}}{\eta}\right)^{d_{T-1}} 
        }{\delta_T}\Bigg)
        }{3NC_ph^p\underline{f}\underline{k}}+\sqrt{ \frac{2\alpha\overline{k}
        \log \Bigg(\frac{\mc
         O(1) \left(\frac{D_{T-1}}{\eta}\right)^{d_{T-1}} 
        }{\delta_T}\Bigg)
        \left(
        \sigma^2_T+3M_Th
        \right)}{NC_ph^p\underline{f}\underline{k}} } + 2 L_{T} \eta
\end{IEEEeqnarray*}
with probability at least $1-\delta_T$. Applying union bound over the $N$ sample points in $\{\widetilde{\bm \gamma}_{i}^{(a)}\}^N_{i=1}$, for any fixed $\x'\in\mathcal{X}_{T-2}(\cdot)$ we have 
 \begin{equation}\label{eq:out_of_sample_T_main}
     \begin{aligned}
           & \left|Q_{T-1}(\bm x',\bm \xi')  -\widehat Q_{T-1}(\bm x',\bm \xi') \right| 
            \\
      \leq  &  M_T h + \frac{2\alpha
            \log \Bigg(\frac{
             \mc O(1) N \left(\frac{D_{T-1}}{\eta}\right)^{d_{T-1}} 
            }{\delta_T}\Bigg)
            }{3NC_p\underline{f}h^{^p}}+\sqrt{ \frac{2\alpha\overline{k}
            \log \Bigg
            (\frac{
             \mc O(1) N \left(\frac{D_{T-1}}{\eta}\right)^{d_{T-1}} 
            }{\delta_T}
            \Bigg)
            \left(
            \sigma^2_T+3M_Th
            \right)}{NC_p\underline{f}h^{^p}} } + 2 L_{T} \eta,
            \\
            & \qquad\qquad\qquad\qquad\qquad\qquad
        \qquad\qquad\qquad\qquad\qquad\qquad
        \qquad\qquad\qquad\qquad\quad\;\;
        \forall \bm\xi'\in \{\widetilde{\bm \gamma}_i^{(a)}\}^N_{i=1}
    \end{aligned}
\end{equation}
with probability at least $1-\delta_T$.

\noindent \underline{Stage $t = T-1$}: From \Cref{corol:bound_for_a_fixed}, for any fixed $\x\in\mathcal{X}_{T-2}(\cdot)$, $\x'\in\mathcal{X}_{T-3}(\cdot)$, and $\bm\xi'\in\Xi$, we have
\begin{multline}
    \left| \EE_{\bm \xi'} \left[ Q_{T-1}(\bm x,\xit) \right] - \widehat \EE_{\bm \xi'} \left[  Q_{T-1}(\bm x,\xit)\right] \right| 
    \\ 
        \leq M_{T-1} h + \frac{2\alpha\log\left(\frac{2}{\delta_{T-1}}\right)}{3NC_ph^p\underline{f}\underline{k}}+\sqrt{ \frac{2\alpha\overline{k}\log\left(\frac{2}{\delta_{T-1}}\right)\left(
        \sigma^2_{T-1}
        +3M_{T-1}h
        \right)}{NC_ph^p\underline{f}\underline{k}} }
        \label{eq:error_bound_T_1_main}
\end{multline}
with probability at least $1-\delta_{T-1}$.
Note that $ \widehat{\EE}_{\bm{\xi}'} [ Q_{T-1}(\cdot)] $ in \eqref{eq:error_bound_T_1_main} represents the approximate expectation of the \emph{true} value function $ Q_{T-1}(\cdot) $, which cannot be computed without access to the conditional distribution. Instead, in the MRST problem, we need to evaluate $ \widehat{\EE}_{\bm{\xi}'} [ \widehat{Q}_{T-1}(\cdot)]$.
Therefore, replacing $Q_{T-1}(\cdot)$ with $\widehat{Q}_{T-1}(\cdot)$ using $ \eqref{eq:out_of_sample_T_main} $ introduces additional errors. Specifically, by applying the union bound to \eqref{eq:out_of_sample_T_main} and \eqref{eq:error_bound_T_1_main}, for any fixed \(\x \in \mathcal{X}_{T-2}(\cdot)\), \(\x' \in \mathcal{X}_{T-3}(\cdot)\), and \(\bm\xi' \in \Xi\), we obtain
\begin{align}
\label{eq:gg_guarantee}
    &\left| \EE_{\bm \xi'} \left[ Q_{T-1}(\bm x, \xit) \right] - \widehat \EE_{\bm \xi'} \left[ \widehat Q_{T-1}(\bm x, \xit) \right] \right| \nonumber \\
    \leq & M_T h + \frac{2\alpha
        \log \Bigg(\frac{
         \mc O(1) N \left(\frac{D_{T-1}}{\eta}\right)^{d_{T-1}} 
        }{\delta_T}\bigg)
        }{3NC_ph^p\underline{f}\underline{k}}+\sqrt{ \frac{2\alpha\overline{k}
        \log \Bigg(\frac{\mc
         O(1) N \left(\frac{D_{T-1}}{\eta}\right)^{d_{T-1}}
        }{\delta_T}\Bigg)
        \left(
        \sigma^2_T+3M_Th
        \right)}{NC_ph^p\underline{f}\underline{k}} } + 2 L_{T} \eta \nonumber \\
       &+M_{T-1} h + \frac{2\alpha\log\left(\frac{2}{\delta_{T-1}}\right)}{3NC_ph^p\underline{f}\underline{k}}+\sqrt{ \frac{2\alpha\overline{k}\log\left(\frac{2}{\delta_{T-1}}\right)\left(
        \sigma^2_{T-1}
        +3M_{T-1}h
        \right)}{NC_ph^p\underline{f}\underline{k}} }
\end{align}

\noindent with probability at least $1-\delta_T-\delta_{T-1}$. Similar to \Cref{thm:Generalization Bound_for_a_Continuous}, we construct an $\eta-$net of the feasible state space, i.e., a set of finite points $\mathcal{X}_{T-2}^{\eta}(\cdot)  \subset \mathcal{X}_{T-2} (\cdot) $ such that for any $\bs x \in \mathcal{X}_{T-2} (\cdot)$, there exists $\bs y \in \mathcal{X}_{T-2}^{\eta} (\cdot)$ such that $\Vert \bs x - \bs y \Vert \leq \eta$. 
Given a finite diameter $D_{T-2}$ as stated in \ref{a3}, the cardinality of the set $\vert \mathcal{X}_{T-2}^{\eta}(\cdot) \vert = \mc O(1)(D_{T-2}/\eta)^{d_{T-2}}$ is controlled by $\eta > 0$.
Since the value function $Q_{T-1}(\cdot)$ is $L_{T-1}$-Lipschitz continuous in $\bs x$, \Cref{corol:expect_lipschitz} implies that for any $\bm x \in \mc X_{T-2}(\cdot)$, there exists $\bs y \in \mathcal{X}_{T-2}^{\eta} (\cdot )$, $\Vert \bs x - \bs y \Vert \leq \eta$, such that:
\begin{equation}
\label{eq:lipschitz_T}
\left| \EE_{\bm \xi'} \left[ Q_{T-1}(\bm x, \xit)  \right] - \EE_{\bm \xi'} \left[ Q_{T-1}(\bm y, \xit)\right] \right| \leq L_{T-1}\eta.
\end{equation}
Furthermore, applying union bound to \eqref{eq:gg_guarantee} over all $\bm y\in \mc X_{T-2}^{\eta}(\cdot)$, for fixed $\bm x'\in \mc X_{T-3}(\cdot)$ and $\bm\xi'\in\Xi$ we get 
\begin{align*}
    &\left| \EE_{\bm \xi'} \left[ Q_{T-1}(\bm y, \xit)  \right] - \widehat \EE_{\bm \xi'} \left[ \widehat Q_{T-1}(\bm y, \xit)\right] \right| \\
     \leq& M_T h + \frac{2\alpha
        \log \Bigg(\frac{
         \mc O(1) N \left(\frac{D_{T-1}}{\eta}\right)^{d_{T-1}}
         \left(\frac{D_{T-2}}{\eta}\right)^{d_{T-2}}
        }{\delta_T}\Bigg)
        }{3NC_ph^p\underline{f}\underline{k}}\\
        &+
        \sqrt{ \frac{2\alpha\overline{k}
        \log \Bigg(\frac{
         \mc O(1) N \left(\frac{D_{T-1}}{\eta}\right)^{d_{T-1}}
         \left(\frac{D_{T-2}}{\eta}\right)^{d_{T-2}} 
        }{\delta_T}\Bigg)
        \left(
        \sigma^2_T+3M_Th
        \right)}{NC_ph^p\underline{f}\underline{k}} } + 2 L_{T} \eta + M_{T-1} h
        \\
        & + \frac{2\alpha\log\Bigg(\frac{
         \mc O(1) 
         \left(\frac{D_{T-2}}{\eta}\right)^{d_{T-2}}
        }{\delta_{T-1}}\Bigg)}{3NC_ph^p\underline{f}\underline{k}} + \sqrt{ \frac{2\alpha\overline{k}\log\Bigg(\frac{
         \mc O(1) 
         \left(\frac{D_{T-2}}{\eta}\right)^{d_{T-2}}
        }{\delta_{T-1}}\Bigg)\left(
        \sigma^2_{T-1}
        +3M_{T-1}h
        \right)}{NC_ph^p\underline{f}\underline{k}} },
        \qquad 
        \\
        & \qquad\qquad\qquad\qquad\qquad\qquad
        \qquad\qquad\qquad\qquad\qquad\qquad
        \qquad\qquad\;\;\;\;
        \forall \bs y \in \mathcal{X}_{T-2}^{\eta} \left( \bm x', \bm \xi' \right) 
\end{align*}
Using the Lipschitz continuity of $Q_{T-1}(\cdot)$ and $\widehat Q_{T-1}(\cdot)$, we obtain
\begin{align*}
    &\left| \EE_{\bm \xi'} \left[ Q_{T-1}(\bm x, \xit)  \right] - \widehat \EE_{\bm \xi'} \left[ \widehat Q_{T-1}(\bm x, \xit)\right] \right| \\
     \leq&M_T h + \frac{2\alpha
        \log \Bigg(\frac{
         \mc O(1) N \left(\frac{D_{T-1}}{\eta}\right)^{d_{T-1}}
         \left(\frac{D_{T-2}}{\eta}\right)^{d_{T-2}}
        }{\delta_T}\Bigg)
        }{3NC_ph^p\underline{f}\underline{k}}\\
        &+
        \sqrt{ \frac{2\alpha\overline{k}
        \log \Bigg(\frac{
         \mc O(1) N \left(\frac{D_{T-1}}{\eta}\right)^{d_{T-1}}
         \left(\frac{D_{T-2}}{\eta}\right)^{d_{T-2}} 
        }{\delta_T}\Bigg)
        \left(
        \sigma^2_T+3M_Th
        \right)}{NC_ph^p\underline{f}\underline{k}} } + 2 L_{T} \eta + M_{T-1} h
        \\
        & + \frac{2\alpha\log\Bigg(\frac{
         \mc O(1) 
         \left(\frac{D_{T-2}}{\eta}\right)^{d_{T-2}}
        }{\delta_{T-1}}\Bigg)}{3NC_ph^p\underline{f}\underline{k}} + \sqrt{ \frac{2\alpha\overline{k}\log\Bigg(\frac{
         \mc O(1) 
         \left(\frac{D_{T-2}}{\eta}\right)^{d_{T-2}}
        }{\delta_{T-1}}\Bigg)\left(
        \sigma^2_{T-1}
        +3M_{T-1}h
        \right)}{NC_ph^p\underline{f}\underline{k}} }+2L_{T-1}\eta,\\
        &
        \qquad\qquad\qquad\qquad\qquad\qquad
        \qquad\qquad\qquad\qquad\qquad\qquad
        \qquad\qquad\qquad\qquad\quad\;\;
        \forall \bs x \in \mathcal{X}_{T-2} \left( \bm x',\bm \xi'\right)         
\end{align*}
with probability at least $1-\delta_T-\delta_{T-1}$.
Then, minimizing over $\mathcal{X}_{T-2} ( \bm x',\bm \xi') $, for any fixed $\bm x' \in \mc X_{T-3}(\cdot)$ and $\bm \xi'\in\Xi$, we have
\begin{IEEEeqnarray*}{lCl}
  &&\bigg| \min_{ \x \in \mc X_{T-2}\left(\bm x', \bm \xi' \right) } \left\lbrace 
  c_{_{T-2}}(\bm x,\bm\xi')
  +  
  \EE_{\bm \xi'} \left[ Q_{T-1}(\bm x, \bm \xit)\right]  \right\rbrace  
  -
  \min_{ \x \in \mc X_{T-2}\left(\bm x', \bm \xi' \right) } \left\lbrace 
  c_{_{T-2}}(\bm x,\bm\xi')
  + \widehat \EE_{\bm \xi'} \left[ \widehat Q_{T-1}(\bm x, \xit)\right]  \right\rbrace \bigg| \\
 &=& \left|Q_{T-2}(\bm x',\bm \xi') - \widehat Q_{T-2}(\bm x',\bm \xi') \right| \\
 & \leq
     & M_T h + \frac{2\alpha
        \log \Bigg(\frac{
         \mc O(1) N \left(\frac{D_{T-1}}{\eta}\right)^{d_{T-1}}
         \left(\frac{D_{T-2}}{\eta}\right)^{d_{T-2}}
        }{\delta_T}\Bigg)
        }{3NC_ph^p\underline{f}\underline{k}}
        \\
        &&+\sqrt{ \frac{2\alpha \overline{k}
        \log \Bigg(\frac{\mc
         O(1) N \left(\frac{D_{T-1}}{\eta}\right)^{d_{T-1}}
         \left(\frac{D_{T-2}}{\eta}\right)^{d_{T-2}} 
        }{\delta_T}\Bigg)
        \left(
        \sigma^2_T+3M_Th
        \right)}{NC_ph^p\underline{f}\underline{k}} } + 2 L_{T} \eta+M_{T-1} h 
        \\
        && + \frac{2\alpha\log\Bigg(\frac{\mc
         O(1) 
         \left(\frac{D_{T-2}}{\eta}\right)^{d_{T-2}}
        }{\delta_{T-1}}\Bigg)}{3NC_ph^p\underline{f}\underline{k}}+\sqrt{ \frac{2\alpha\overline{k}\log\Bigg(\frac{\mc
         O(1) 
         \left(\frac{D_{T-2}}{\eta}\right)^{d_{T-2}}
        }{\delta_{T-1}}\Bigg)\left(
        \sigma^2_{T-1}
        +3M_{T-1}h
        \right)}{NC_ph^p\underline{f}\underline{k}} }+2L_{T-1}\eta
\end{IEEEeqnarray*}
with probability at least $1-\delta_T-\delta_{T-1}$. Finally, applying union bound 
over the $N$ sample points in $\{\widetilde{\bm \gamma}_i^{(b)}\}^N_{i=1}$, for any fixed $\bm x'\in\mc X_{T-3}(\cdot)$, we have
\begin{align*}
    &\left| 
    Q_{T-2}(\bm x', \xi'
    ) 
    - 
    \widehat Q_{T-2}(\bm x', \xi'
    ) 
    \right| \\
     \leq& M_T h + \frac{2\alpha
        \log \Bigg(\frac{\mc
         O(1) N^2 \left(\frac{D_{T-1}}{\eta}\right)^{d_{T-1}}
         \left(\frac{D_{T-2}}{\eta}\right)^{d_{T-2}}
        }{\delta_T}\Bigg)
        }{3NC_ph^p\underline{f}\underline{k}}
        \\
        &+\sqrt{ \frac{2\alpha\overline{k}
        \log \Bigg(\frac{\mc
         O(1) N^2 \left(\frac{D_{T-1}}{\eta}\right)^{d_{T-1}}
         \left(\frac{D_{T-2}}{\eta}\right)^{d_{T-2}} 
        }{\delta_T}\Bigg)
        \left(
        \sigma^2_T+3M_Th
        \right)}{NC_ph^p\underline{f}\underline{k}} } + 2 L_{T} \eta+M_{T-1} h 
        \\
        & +\frac{2\alpha\log\Bigg(\frac{\mc
         O(1) N 
         \left(\frac{D_{T-2}}{\eta}\right)^{d_{T-2}}
        }{\delta_{T-1}}\Bigg)}{3NC_ph^p\underline{f}\underline{k}}
        +
        \sqrt{ \frac{2\alpha\overline{k}\log\Bigg(\frac{\mc
         O(1) N
         \left(\frac{D_{T-2}}{\eta}\right)^{d_{T-2}}
        }{\delta_{T-1}}\Bigg)\left(
        \sigma^2_{T-1}
        +3M_{T-1}h
        \right)}{NC_ph^p\underline{f}\underline{k}} }
        +
        2L_{T-1}\eta, \\
        &
        \qquad\qquad\qquad\qquad\qquad\qquad
        \qquad\qquad\qquad\qquad\qquad\qquad
        \qquad\qquad\qquad\qquad\qquad\quad\;\;
        \forall \bm\xi' \in \{\widetilde{\bm \gamma}_i^{(b)}\}^N_{i=1} 
\end{align*}
with probability at least $1-\delta_T-\delta_{T-1}$.
\noindent Repeating the above steps for \underline{stages $t=T-2$ to $t=2$}, we have the following result: for all $\bm x_1 \in \mathcal{X}_{1} \left( \bm x_0, \bm \xi_1 \right)$, we have
  \begin{multline*}
        \left| 
        \EE_{\bm \xi_1} \left[ Q_2 (\bm x_1, \widetilde{\bm \xi}_{2} )   \right] 
        -
        \widehat \EE_{\bm \xi_1} \left[ \widehat Q_{2}(\bm x_1,  \xit_{2}) \right] 
        \right| \\
        \leq \sum_{t=2}^{T}
        M_th 
        +
        \frac{2\alpha\log \Big(\frac{\mc O(1) N^{t-2} \prod_{s =1}^{t-1} \left(\frac{D_{s}}{\eta}\right)^{d_s} }{ \delta_t }\Big)}{3NC_ph^p\underline{f}\underline{k}}
        +
         \sqrt{ \frac{2\alpha\overline{k}\log \Big(\frac{\mc O(1) N^{t-2} \prod_{s =1}^{t-1} \left(\frac{D_{s}}{\eta}\right)^{d_s} }{ \delta_t }\Big)\left(
        \sigma^2_{t}
        +3M_{t}h
        \right)}{NC_ph^p\underline{f}\underline{k}} }
        +
        2 L_t \eta,
 \end{multline*}
with probability at least $1-\sum_{t=2}^{T} \delta_t$. This completes the proof. 
\end{proof}

\noindent \Cref{thm:out_of_sample} enables us to further derive the suboptimality bound for the MRST problem. In what follows, to provide clear insight, we simplify the bound by taking the largest values of the problem-dependent parameters.
\begin{corollary}
[Suboptimality Bound]
\label{pro:gen_error}
\setlength{\abovedisplayskip}{0pt}%
Define $\displaystyle M = \max_{t\in[T-1]} \; M_{t+1},\;L = \max_{t\in[T-1]} \; L_{t+1},\; D = \max_{t\in[T-1]} D_{t+1}$, $\displaystyle d = \max_{t\in[T-1]} d_{t+1}$, and $\displaystyle \sigma^2 = \max_{t\in [T-1]} \sigma_{t+1}^2$. Suppose that ${\bm x}_1^\star$ is an optimal solution to the MSP problem~\eqref{eq:true_dp_first}, and $\widehat{\bm x}_1^N$ is an optimal solution to the MRST problem~\eqref{eq:approx_dp}, constructed using two independent trajectories, each of length $N+1$. Then, for any fixed $\delta \in (0,1]$, $h\geq0$, and $\eta > 0$, we have
\begin{multline}
        \left|
        \left(
        c_{{1}}(\bm x_1^\star,\bm\xi_1)
        + \EE_{\bm \xi_1} \left[ Q_2 (\bm x_1^\star, \widetilde{\bm \xi}_{2} )  \right] \right) 
        -
        \left( 
        c_{{1}}(\widehat{\bm x}_1^N,\bm\xi_1)
        + \EE_{\bm \xi_1} \left[ Q_{2}(\widehat{\bm x}_1^N,  \xit_{2}) \right] \right) \right| 
         \\  
 \leq
 2T
 \scalebox{1.5}{$\Bigg($}
 Mh 
 + 
\frac{2\alpha T\log\Big( \frac{\mc O(1)NT\left({D}/{\eta}\right)^{d}}{\delta} \Big)}{3NC_ph^p\underline{f}\underline{k}}
 +
 \sqrt{\frac{2\alpha\overline{k} T\log\Big( \frac{\mc O(1)NT\left({D}/{\eta}\right)^{d}}{\delta} \Big)\left( \sigma^2 + 3Mh \right)}{NC_ph^p\underline{f}\underline{k}} } +L\eta
 \scalebox{1.5}{$\Bigg)$}
 \label{eq:sample_comp}
\end{multline}
with probability at least $1-\delta$.
\label{cor_simplified_gen_bound}
\end{corollary}
\begin{proof}
    Using the largest problem-dependent parameters $M$, $L$, $D$, $d$, and $\sigma^2$, we can derive the following upper bound on \eqref{eq_OOS}: for all $\bm x_1 \in \mathcal{X}_{1} \left( \bm x_0, \bm \xi_1 \right)$, we have
    \begin{multline}
        \left|  \left( 
        c_{{1}}(\bm x_1,\bm\xi_1)
        + \EE_{\bm \xi_1} \left[ Q_2 (\bm x_1, \widetilde{\bm \xi}_{2} )   \right] \right) 
        -
        \left( 
        c_{{1}}(\bm x_1,\bm\xi_1)
        + \widehat \EE_{\bm \xi_1} \left[ \widehat Q_{2}(\bm x_1,  \xit_{2})\right] \right)\right| 
         \\  
 \leq
 \underbrace{
 T
 \scalebox{1.5}{$\Bigg($}
 Mh 
 + 
\frac{2\alpha T\log\Big( \frac{\mc O(1)NT\left(\frac{D}{\eta}\right)^{d}}{\delta} \Big)}{3NC_ph^p\underline{f}\underline{k}}
 +
 \sqrt{\frac{2\alpha\overline{k} T\log\Big( \frac{\mc O(1)NT\left(\frac{D}{\eta}\right)^{d}}{\delta} \Big)\left( \sigma^2 + 3Mh \right)}{NC_ph^p\underline{f}\underline{k}}} +L\eta
 \scalebox{1.5}{$\Bigg)$}
 }_{\displaystyle =\epsilon(\delta)}
 \label{eq_sample_comp}
\end{multline}
with probability at least $1-\delta$. As shown above, $\epsilon(\delta)$ is denoted as the right-hand side of inequality~\eqref{eq_sample_comp}.
Then, for $\bm x_1=\widehat{\bm x}_1^N$, with probability $1-\delta$, we have 
\begin{align}
        c_{{1}}(\widehat{\bm x}_1^N,\bm\xi_1)
        + \EE_{\bm \xi_1} \left[ Q_2 (\widehat{\bm x}_1^N, \widetilde{\bm \xi}_{2} )  \right] 
        &\leq
        c_{{1}}(\widehat{\bm x}_1^N,\bm\xi_1)
        + \widehat{\EE}_{\bm \xi_1} \left[ \widehat{Q}_{2}(\widehat{\bm x}_1^N,  \xit_{2}) \right]  +\epsilon(\delta)
        \nonumber
        \\
        &\leq
         c_{{1}}({\bm x}_1^\star,\bm\xi_1)
        + \widehat{\EE}_{\bm \xi_1} \left[ \widehat{Q}_{2}({\bm x}_1^\star,  \xit_{2}) \right]  +\epsilon(\delta).
        \label{eq_idk}
\end{align}
The second inequality \eqref{eq_idk} holds because ${\bm x}_1^\star$ is suboptimal for the MRST problem \eqref{eq:approx_dp}. Likewise, for $\bm x_1 = {\bm x}_1^\star$, with probability $1-\delta$, we have the following:
\begin{equation}\label{eq_idk2}
        -\left(
        c_{{1}}({\bm x}_1^\star,\bm\xi_1)
        + \EE_{\bm \xi_1} \left[ Q_2 ({\bm x}_1^\star, \widetilde{\bm \xi}_{2} ) \right] 
        \right)
        +
        \left(
        c_{{1}}({\bm x}_1^\star,\bm\xi_1)
        + \widehat{\EE}_{\bm \xi_1} \left[ \widehat{Q}_2 ({\bm x}_1^\star, \widetilde{\bm \xi}_{2} )  \right] 
        \right)
        \leq \epsilon(\delta)
        .
\end{equation}
Applying union bound to \eqref{eq_idk} and \eqref{eq_idk2}, we establish \eqref{eq:sample_comp}.
\end{proof}
\noindent 
Note that the bandwidth parameter $h$ in \eqref{eq:sample_comp} should be adjusted based on problem-dependent parameters, including the number of samples $N$ and the time horizon $T$, to obtain a sharper bound.
Since, the inequality \eqref{eq:sample_comp} holds with probability at least $1-\delta$ for any $h \geq 0$, we can consider the right-hand side of \eqref{eq:sample_comp} as a minimization problem over $h$ to obtain a sharper bound. 
While applying Bernstein's inequality in \eqref{eq:bernstein_placeholder} generally provides a tighter suboptimality bound compared to Hoeffding's inequality~(\citet[Theorem 2.2.5]{vershynin2018high}), it leads to a more complex optimization problem. 
In contrast, the Hoeffding version of the bound leads to a simpler optimization problem, allowing us to derive an analytical formula for the optimal bandwidth $h^\star$, as follows:
\begin{align}
     h^\star
     &=\left( 
     \frac{p}{2M}\cdot
     \sqrt{
     \frac{\Gamma(p/2+1)}{\pi^{p/2}}\cdot
     \frac{\alpha \overline{k} T}{\underline{f}\underline{k}N}\log\left(\frac{N T(D/\eta)^d}{\delta}\right)
     }
       \right)^{\frac{2}{p+2}}
       \label{eq:optimal_band_pre}
       \\
       &\approx
       \left( 
     \frac{p}{2M}\cdot
     \sqrt{
     \frac{\sqrt{\pi p}\left(\frac{p}{2e}\right)^{{p}/{2}}}{\pi^{{p}/{2}}}\cdot
     \frac{\alpha \overline{k} T}{\underline{f}\underline{k}N}\log\left(\frac{N T(D/\eta)^d}{\delta}\right)
     }
       \right)^{\frac{2}{p+2}}.
    \label{eq:optimal_band}
\end{align}
In \eqref{eq:optimal_band_pre}, we express $1/C_p$ explicitly as $\frac{\Gamma(p/2+1)}{\pi^{p/2}}$ (see \Cref{prop:kernel_estimation_error}), and in \eqref{eq:optimal_band}, the gamma function $\Gamma(p/2+1)$ is approximated using Stirling's formula~(\citet[Equation 1]{robbins1955remark}). We adopt \eqref{eq:optimal_band} to derive a closed-form of the suboptimality bound~\eqref{eq:sample_comp}.
By substituting \eqref{eq:optimal_band} into \eqref{eq:sample_comp}, we can further simplify the suboptimality bound to 
$\widetilde{\mathcal{O}}(T^{\frac{p+3}{p+2}}N^{-\frac{1}{p+2}})$. 
Equivalently, the sample complexity of the MRST problem is
$\widetilde{\mathcal{O}}(T^{p+3} \epsilon^{-p-2})$
in order to obtain an $\epsilon$-optimal solution to the MSP problem \eqref{eq:true_dp_first} as claimed in \Cref{corol:sample_comp}.
\begin{remark}
    We compare our sample complexity 
    $\widetilde{\mathcal{O}}(T^{p+3} \epsilon^{-p-2})$
    with that of the SAA method. 
    We remark that the additional dependence on $p$ arises because the conditional distribution is unknown in our data-driven setup, requiring us to estimate it from data. 
    \citet{shapiro2006complexity} analyzes the SAA method for the MSP problem under a known conditional distribution and demonstrates that its suboptimality bound is $\widetilde{\mathcal{O}}(T^{\frac{1}{2}} N^{-\frac{1}{2T}})$, leading to a sample complexity rate of $\widetilde{\mathcal{O}}(T^T \epsilon^{-2T})$, which suffers from the curse of dimensionality with respect to $T$. Since a known distribution allows them to generate as many samples as needed to improve performance, the scarcity of samples is not their primary concern. Rather, the authors highlight a significant computational limitation of applying the SAA method to the MSP problem: as $T$ increases, the computational and memory demands to obtain an $\epsilon$-optimal solution grow exponentially. This makes SAA impractical for the MSP problem \eqref{eq:true_dp} when the planning horizon is large. The exponentially increasing computational burden of sampling-based methods in solving the MSP problem under Markov uncertainty is a key challenge. 
    Whether in a data-driven setting or not, most existing methods circumvent the curse of dimensionality by imposing more stringent assumptions. 
    For instance, they often assume that stagewise independence in the data process (\citet{pereira1991multi,shapiro2011analysis}), i.e., a random vector $\widetilde{\bm\xi}_{t+1}$ is independent of the history $\bm \xi_{[1:t]}$, ignoring temporal dependence in a data process that naturally arises in many applications. While some methods attempt to address this issue by incorporating Markovian uncertainty, they often restrict it solely to the right-hand side of constraints to ensure the MSP problem remain tractable (\citet{infanger1996cut,shapiro2013risk}). Unfortunately, these assumptions limit the applicability of such methods in addressing various real-world problems with a large number of decision stages.
    We acknowledge the existence of several works addressing general Markovian uncertainty that are based on a discretization method known as optimal quantization (\citet{bonnans2012energy, castro2022markov, lohndorf2019modeling, philpott2012dynamic}). However, there has been no rigorous sample complexity analysis that addresses the curse of dimensionality. 
\end{remark}

\section{Numerical Experiments}
In this section, we conduct stylized numerical experiments to demonstrate the effectiveness of the proposed method. The experiments are implemented using Python 3.7 on a 6-core, 2.3GHz Intel Core i7 CPU laptop with 16GB RAM.
\begin{figure}[htbp]
  \centering
  \includegraphics[width=1.0\textwidth]{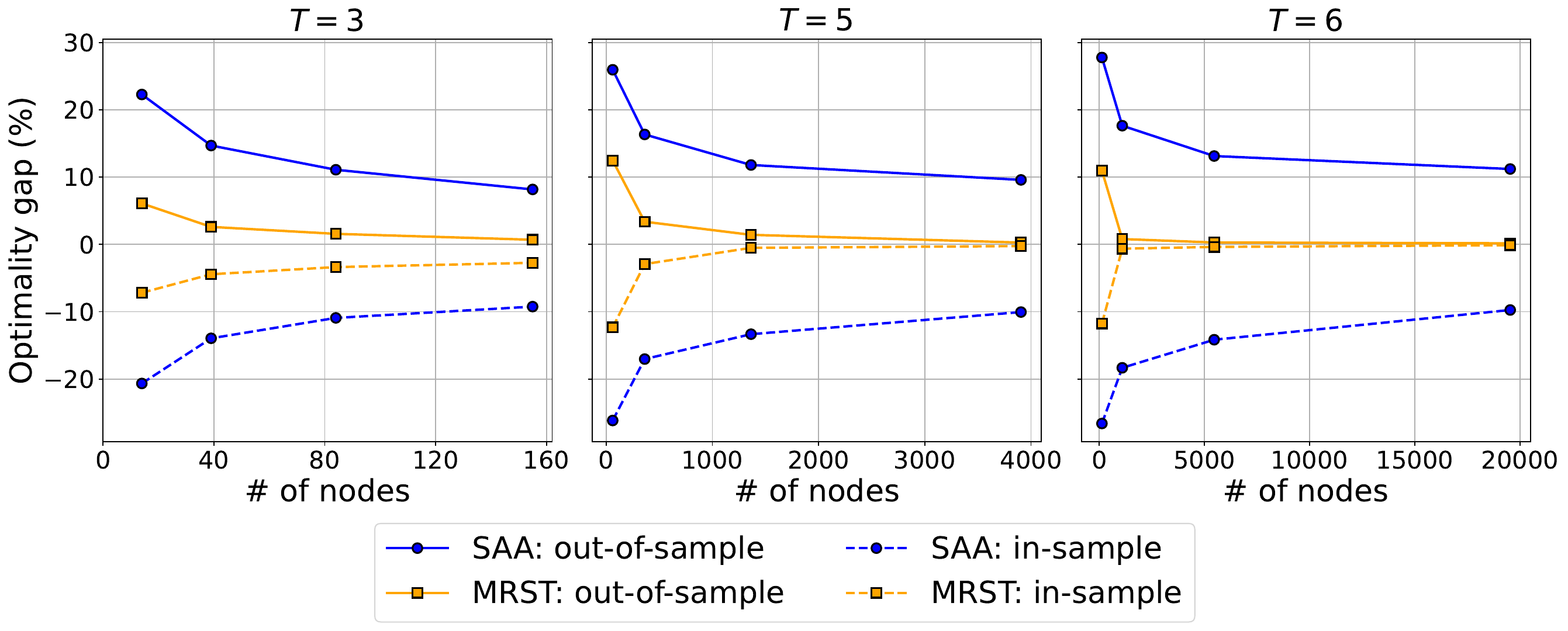}  
    \caption{Comparison of optimality gaps between SAA and MRST: solid lines represent the out-of-sample suboptimality for MRST (blue) and SAA (yellow), while dashed lines represent the in-sample suboptimality.}
\label{fig:suboptimallity}
\end{figure}
We consider the discrete-time Linear Quadratic Gaussian (LQG) control problem, where the objective is to minimize the expected quadratic cost of a controller operating in a linear dynamic environment with additive Gaussian noise. Specifically, the system dynamics are given by 
$$
\bm{x}_{t} = \bm{A} \bm{x}_{t-1} + \bm{B} \bm{u}_{t} + \bm{\xi}_t,
$$
where \(\bm{A} \in \mathbb{R}^{n \times n}\) and \(\bm{B} \in \mathbb{R}^{n \times m}\) describe the linear dynamics. For each stage \(t \in [T-1]\), the cost is a quadratic function of the state $\bm{x}_{t-1} \in \mathbb{R}^{n}$ and control $\bm{u}_{t} \in \mathbb{R}^{m}$, expressed as 
$$
\bm{x}_{t-1}^\top \bm{G} \bm{x}_{t-1}^{} + \bm{u}_{t}^\top \bm{H} \bm{u}_{t}^{},
$$
where $\bm{G} \in \mathbb{R}^{n \times n}$ and $\bm{H} \in \mathbb{R}^{m \times m}$ are weighting matrices for the state and control, respectively. The terminal cost is given by 
$
\bm{x}_{T}^\top \bm{G} \bm{x}_{T}^{}.
$

 In this problem, we can explicitly split the decision vector into state variables $\bm x_{t-1}$ for each $t\in[T+1]$ and control variables $\bm u_{t}$ for each $t\in[T]$. Then, similar to \eqref{eq:true_dp_first} and \eqref{eq:true_dp}, we can define the the value function as follows:
\begin{equation}
\label{eq:true_dp_first_lqg}
\begin{array}{ccl}
Q_1(\bm x_{0}, \bm\xi_1) = & \textup{min} & \displaystyle 
\bm x_0^\top G \bm x_0^{} + \bm u_1^\top H \bm u_1^{} +  \EE_{\bm \xi_1} \Big[ Q_{2}({\bm x}_{1}, \widetilde{\bm \xi}_{2} )  \Big] \\
& \textup{s.t.} & \displaystyle \bm x_1\in\mathbb{R}^{n}, \;\; \bm u_1\in\mathbb{R}^{m},
\\
&  & \displaystyle \bm x_1 = \bm A \bm x_{0} + \bm B \bm u_{1} + \bm \xi_1,
\end{array}
\end{equation}
and, for each subsequent stage $t\in [T]\setminus \{1\}$, $\bm x_{t-1}$, and $\bm\xi_t$, 
\begin{equation}
\label{eq:true_dp_lqg}
\begin{array}{ccl}
Q_t(\bm x_{t-1},\bm \xi_{t}) = & \textup{min} & \displaystyle 
 \bm x_{t-1}^\top G \bm x_{t-1}^{} + \bm u_t^\top H \bm u_t^{} +  \EE_{\bm \xi_t} \Big[ Q_{t+1}({\bm x}_{t} ,\widetilde{\bm \xi}_{t+1} ) \Big] \\
& \textup{s.t.} & \displaystyle \bm x_t\in\mathbb{R}^{n}, \;\; \bm u_t\in\mathbb{R}^{m},
\\&& \displaystyle \bm x_t = \bm A \bm x_{t-1} + \bm B \bm u_{t} + \bm \xi_{t+1}.
\end{array}
\end{equation}
A key distinction from the standard LQG problem~(\citet{shaiju2008formulas}), adapted to our setting, is that we assume the noise process follows a Markov structure, given by
$$
\bm{\xi}_{t+1} = \bm{C} \bm{\xi}_{t} + \bm{\eta}_{t+1},
$$ 
where $\bm{C} \in \mathbb{R}^{n \times n}$ and $\bm{\eta}_{t+1} \sim \mathcal{N}(\bm{0}, \bm{\Sigma})$. We adopt this problem because, similar to the standard LQG framework, the dynamic programming equations \eqref{eq:true_dp_first_lqg} and \eqref{eq:true_dp_lqg} are quadratic programs and yield optimal solutions and value functions in closed form based on the matrices $\bm A$, $\bm B$, $\bm G$, $\bm H$, $\bm C$, and $\bm \Sigma$, that govern the dynamics and cost. This enables us to compute the \emph{exact} suboptimality of the optimal solutions to the SAA and MRST methods, allowing a comparison of their performance across varying sample sizes and time horizons. 

We conducted 1,000 independent in-sample and out-of-sample tests for SAA and MRST across varying time horizons $T \in \{3,5,6\}$. In \Cref{fig:suboptimallity}, we plot the mean optimality gaps of the in-sample and out-of-sample costs of the two methods, defined as follows:
{\small
\begin{align}
    &\textbf{OptGap}_{\text{in-sample}} = \frac{\text{Mean In-sample Objective Value} - \text{True Optimal Objective Value} }{\text{True Optimal Objective Value}} \times 100, \nonumber \\
    \nonumber  \\
    &\textbf{OptGap}_{\text{out-of-sample}} = \frac{\text{Mean Out-of-sample Objective Value} - \text{True Optimal Objective Value} }{\text{True Optimal Objective Value}} \times 100. \nonumber
\end{align}
}
The x-axis in \Cref{fig:suboptimallity} represents the number of nodes in the scenario tree, i.e., the sample complexity. As expected, the results show that SAA suffers from the curse of dimensionality: to maintain a suboptimality gap within $\pm10\%$, the size of the SAA problem grows exponentially as the time horizon increases. In contrast, MRST achieves a tighter suboptimality gap than SAA, with the gap tighter than $\pm1\%$ when the number of nodes exceeds 4,000. The improvement in the suboptimality gap over an increasing time horizon demonstrates the polynomial sample complexity of our method.

Additionally, obtaining an implementable SAA policy in the out-of-sample test requires re-solving a new instance of the SAA problem at each time step, which demands significant computational effort. For instance, when $T=6$ and the number of conditional samples $N=5$, the number of generated samples is approximately 20,000. 
Our experiments not only align with the literature in demonstrating that the SAA method is impractical for the MSP problem but also show that our method provides a more efficient data-driven alternative to SAA, effectively addressing the curse of dimensionality.

\begin{remark}
For the LQG problem, the exact optimal solution is readily available, hence, using sampling-based approaches does not make sense in practice. Our objective in the numerical experiments is solely to evaluate and compare the suboptimality of the SAA and MRST methods. For a general MSP problem \eqref{eq:true_nested} where closed-form solutions are not available, one may need to adopt Stochastic Dual Dynamic Programming (SDDP) (\citet{pereira1991multi}) and its variations~(\citet{ahmed2003multi,duque2020distributionally,huang2017study,park2022data,philpott2018distributionally,silva2021data,zhang2022stochastic,zou2019stochastic}) to solve the approximate problems. In fact, SDDP is a solution framework compatible with the MRST method. We do not cover this approach, as it is beyond the scope of this paper.
\end{remark}

\section{Conclusion}
In this work, we introduced a novel data-driven approximation method for the MSP problem under a continuous-state Markov process. Compared to the SAA method, our approach not only allows for the implementable policy but also accommodates a broader class of MSP problems with computational efficiency.

Most notably, our method exhibits mild sample complexity with respect to the time horizon, which effectively addresses the curse of dimensionality inherent in the SAA method for solving the MSP problem under Markovian uncertainty. Numerical experiments on the LQG control problem corroborate our theoretical results, demonstrating significant improvements over the SAA method across varying time horizons. These results highlight the potential of our method as a more efficient data-driven alternative to SAA for the MSP problem.

Future research could explore combining our approach with the SDDP algorithm as a solution scheme and investigate real-world applications where data collection is inherently limited. Such efforts would further enhance the practical applicability of our method.

\section*{Acknowledgments}
This work is supported by the National Science Foundation grants 2342505, 2343869, and 2404413.

\bibliographystyle{apalike}
\bibliography{references}

\begin{thebibliography}{}

\bibitem[Ahmed et~al., 2003]{ahmed2003multi}
Ahmed, S., King, A.~J., and Parija, G. (2003).
\newblock A multi-stage stochastic integer programming approach for capacity expansion under uncertainty.
\newblock {\em Journal of Global Optimization}, 26:3--24.

\bibitem[Bellman, 1957]{bellman1957dynamic}
Bellman, R. (1957).
\newblock {\em Dynamic Programming}.
\newblock Princeton University Press.

\bibitem[Bonnans et~al., 2012]{bonnans2012energy}
Bonnans, J.~F., Cen, Z., and Christel, T. (2012).
\newblock Energy contracts management by stochastic programming techniques.
\newblock {\em Annals of Operations Research}, 200:199--222.

\bibitem[Castro et~al., 2022]{castro2022markov}
Castro, M.~P., Bodur, M., and Song, Y. (2022).
\newblock Markov chain-based policies for multi-stage stochastic integer linear programming with an application to disaster relief logistics.
\newblock {\em arXiv preprint arXiv:2207.14779}.

\bibitem[Duque and Morton, 2020]{duque2020distributionally}
Duque, D. and Morton, D.~P. (2020).
\newblock Distributionally robust stochastic dual dynamic programming.
\newblock {\em SIAM Journal on Optimization}, 30(4):2841--2865.

\bibitem[Hanasusanto et~al., 2016]{hanasusanto2016comment}
Hanasusanto, G.~A., Kuhn, D., and Wiesemann, W. (2016).
\newblock A comment on “computational complexity of stochastic programming problems”.
\newblock {\em Mathematical Programming}, 159:557--569.

\bibitem[Huang et~al., 2017]{huang2017study}
Huang, J., Zhou, K., and Guan, Y. (2017).
\newblock A study of distributionally robust multistage stochastic optimization.
\newblock {\em arXiv preprint arXiv:1708.07930}.

\bibitem[Infanger and Morton, 1996]{infanger1996cut}
Infanger, G. and Morton, D.~P. (1996).
\newblock Cut sharing for multistage stochastic linear programs with interstage dependency.
\newblock {\em Mathematical Programming}, 75(2):241--256.

\bibitem[Jiang and Li, 2021]{jiang2021complexity}
Jiang, J. and Li, S. (2021).
\newblock On complexity of multistage stochastic programs under heavy tailed distributions.
\newblock {\em Operations Research Letters}, 49(2):265--269.

\bibitem[Kallenberg, 1997]{kallenberg1997foundations}
Kallenberg, O. (1997).
\newblock {\em Foundations of modern probability}, volume~2.
\newblock Springer.

\bibitem[Kleywegt et~al., 2002]{kleywegt2002sample}
Kleywegt, A.~J., Shapiro, A., and Homem-de Mello, T. (2002).
\newblock The sample average approximation method for stochastic discrete optimization.
\newblock {\em SIAM Journal on optimization}, 12(2):479--502.

\bibitem[L{\"o}hndorf and Shapiro, 2019]{lohndorf2019modeling}
L{\"o}hndorf, N. and Shapiro, A. (2019).
\newblock Modeling time-dependent randomness in stochastic dual dynamic programming.
\newblock {\em European Journal of Operational Research}, 273(2):650--661.

\bibitem[Lozier, 2003]{lozier2003nist}
Lozier, D.~W. (2003).
\newblock Nist digital library of mathematical functions.
\newblock {\em Annals of Mathematics and Artificial Intelligence}, 38:105--119.

\bibitem[Park et~al., 2022]{park2022data}
Park, H., Jia, Z., and Hanasusanto, G.~A. (2022).
\newblock Data-driven stochastic dual dynamic programming: Performance guarantees and regularization schemes.
\newblock {\em Available at Optimization Online}.

\bibitem[Pereira and Pinto, 1991]{pereira1991multi}
Pereira, M.~V. and Pinto, L.~M. (1991).
\newblock Multi-stage stochastic optimization applied to energy planning.
\newblock {\em Mathematical Programming}, 52(1):359--375.

\bibitem[Philpott and De~Matos, 2012]{philpott2012dynamic}
Philpott, A.~B. and De~Matos, V.~L. (2012).
\newblock Dynamic sampling algorithms for multi-stage stochastic programs with risk aversion.
\newblock {\em European Journal of operational research}, 218(2):470--483.

\bibitem[Philpott et~al., 2018]{philpott2018distributionally}
Philpott, A.~B., de~Matos, V.~L., and Kapelevich, L. (2018).
\newblock Distributionally robust sddp.
\newblock {\em Computational Management Science}, 15:431--454.

\bibitem[Reaiche, 2016]{reaiche2016note}
Reaiche, M. (2016).
\newblock A note on sample complexity of multistage stochastic programs.
\newblock {\em Operations Research Letters}, 44(4):430--435.

\bibitem[Robbins, 1955]{robbins1955remark}
Robbins, H. (1955).
\newblock A remark on stirling's formula.
\newblock {\em The American mathematical monthly}, 62(1):26--29.

\bibitem[Shaiju and Petersen, 2008]{shaiju2008formulas}
Shaiju, A. and Petersen, I.~R. (2008).
\newblock Formulas for discrete time lqr, lqg, leqg and minimax lqg optimal control problems.
\newblock {\em IFAC Proceedings Volumes}, 41(2):8773--8778.

\bibitem[Shapiro, 2006]{shapiro2006complexity}
Shapiro, A. (2006).
\newblock On complexity of multistage stochastic programs.
\newblock {\em Operations Research Letters}, 34(1):1--8.

\bibitem[Shapiro, 2011]{shapiro2011analysis}
Shapiro, A. (2011).
\newblock Analysis of stochastic dual dynamic programming method.
\newblock {\em European Journal of Operational Research}, 209(1):63--72.

\bibitem[Shapiro et~al., 2021]{shapiro2021lectures}
Shapiro, A., Dentcheva, D., and Ruszczynski, A. (2021).
\newblock {\em Lectures on stochastic programming: modeling and theory}.
\newblock SIAM.

\bibitem[Shapiro and Nemirovski, 2005]{shapiro2005complexity}
Shapiro, A. and Nemirovski, A. (2005).
\newblock On complexity of stochastic programming problems.
\newblock {\em Continuous optimization: Current trends and modern applications}, pages 111--146.

\bibitem[Shapiro et~al., 2013]{shapiro2013risk}
Shapiro, A., Tekaya, W., da~Costa, J.~P., and Soares, M.~P. (2013).
\newblock Risk neutral and risk averse stochastic dual dynamic programming method.
\newblock {\em European journal of operational research}, 224(2):375--391.

\bibitem[Silva et~al., 2021]{silva2021data}
Silva, T., Vallad{\~a}o, D., and Homem-de Mello, T. (2021).
\newblock A data-driven approach for a class of stochastic dynamic optimization problems.
\newblock {\em Computational Optimization and Applications}, 80:687--729.

\bibitem[Tibshirani and Wasserman, 2013]{tibshirani2013nonparametric}
Tibshirani, R. and Wasserman, L. (2013).
\newblock Nonparametric regression.
\newblock {\em Statistical Machine Learning, Spring}.

\bibitem[Vershynin, 2018]{vershynin2018high}
Vershynin, R. (2018).
\newblock {\em High-dimensional probability: An introduction with applications in data science}, volume~47.
\newblock Cambridge university press.

\bibitem[Zhang and Sun, 2022]{zhang2022stochastic}
Zhang, S. and Sun, X.~A. (2022).
\newblock Stochastic dual dynamic programming for multistage stochastic mixed-integer nonlinear optimization.
\newblock {\em Mathematical Programming}, 196(1):935--985.

\bibitem[Zou et~al., 2019]{zou2019stochastic}
Zou, J., Ahmed, S., and Sun, X.~A. (2019).
\newblock Stochastic dual dynamic integer programming.
\newblock {\em Mathematical Programming}, 175:461--502.

\end{thebibliography}

\end{document}